\documentclass[11pt]{amsart}

\usepackage{graphicx}

\usepackage{subcaption} 
\usepackage{xfrac}   
\usepackage{faktor}
\usepackage{tikz-cd}
\usepackage{mathrsfs}
\usepackage{hyperref}
\usepackage{amsmath}
\usepackage{amssymb}
\hypersetup{bookmarks=true,
	unicode=true,
	colorlinks=true,
	citecolor=black,
	linkcolor=black,
	urlcolor=black,
	plainpages=false,
	pdfpagelabels=true}

 \usepackage{url}	
 \allowdisplaybreaks 

\usepackage{tikz-cd}
\usepackage{pgf}

\usepackage{xcolor}

\usepackage{comment} 

\usepackage[all]{xy}
\usepackage{enumitem}

\newtheorem{teo}{Theorem}[section]
\newtheorem{theorem}[teo]{Theorem}

\newtheorem{corollary}[teo]{Corollary}

\newtheorem{lemma}[teo]{Lemma}

\newtheorem{proposition}[teo]{Proposition}

\theoremstyle{definition}
\newtheorem{definition}[teo]{Definition}
\newtheorem{defi}[teo]{Definition}
\newtheorem{example}[teo]{Example}

\newtheorem*{mainthm}{Main Results}

\theoremstyle{remark}
\newtheorem{remark}[teo]{Remark}

\numberwithin{figure}{section}

\newcommand{\CC}{\mathcal{C}}

\newcommand{\DD}{\mathcal{D}}

\DeclareMathOperator{\Mon}{Mon}

\DeclareMathOperator{\hocolim}{hocolim}

\DeclareMathOperator{\hocofib}{hocofib}

\DeclareMathOperator{\dom}{dom}

\DeclareMathOperator{\cod}{cod}

\DeclareMathOperator{\sd}{sd}

\newcommand{\C}{\mathcal{C}}

\newcommand{\T}{\mathcal{T}}

\newcommand{\Ran}{\mathrm{Ran}}

\newcommand{\Lan}{\mathrm{Lan}}

\newcommand{\Fun}{\mathrm{Fun}}

\DeclareMathOperator*{\colim}{colim}

\newcommand{\Ob}{\mathrm{Ob}}

\newcommand{\Set}{\mathbf{Set}}

\newcommand{\Cat}{\mathbf{Cat}}

\newcommand{\Hot}{\mathbf{Hot}}

\newcommand{\Dinj}{\Delta_{\mathrm{inj}}}
\newcommand{\Hom}{\mathrm{Hom}}

\newcommand{\op}{\mathrm{op}}

\newcommand{\Top}{\mathrm{Top}}

\newcommand{\A}{\mathcal{A}}

\newcommand{\B}{\mathrm{B}}

\newcommand{\EE}{\mathcal{E}}

\newcommand{\FF}{\mathcal{F}}

\begin{document}

\title[One monoid to realize them all]{One monoid to realize them all \\ A master monoid for homotopy types}
\thanks{}

\author[Isaac Carcacía-Campos]{%
	Isaac Carcacía-Campos%
}

 \address{%
           	Isaac Carcacía-Campos
            \\
              Departamento de Matem\'aticas, Universidade de Santiago de Compostela, 15782-SPAIN}
               \email{isaac.c.campos@usc.es}

\begin{abstract}

Homotopy theory has traditionally been approached either through topological spaces or through categories of presheaves on suitable combinatorial categories, a perspective subsumed by Grothendieck's theory of test categories. In this paper, we develop a more algebraic viewpoint: homotopy types can be represented by actions of a single monoid on sets.

This monoid is obtained by identifying all the objects of the semisimplicial category while retaining its face operators. Every CW homotopy type is the classifying space of the translation category of one of its actions, and its truncation at dimension $n$ is finitely presented and realizes exactly the homotopy types of dimension at most $n$.

The construction rests on a general result: every set-valued functor on a small category freely generates an action of its universal monoid, and its category of elements is a full coreflective subcategory of the resulting translation category. Hence every weak test category gives rise to a monoid whose actions realize all homotopy types; when the category has more than one object, however, this monoid is not itself a weak test category, so that a single monoid may realize all homotopy types without modelling homotopy theory as a localization. The obstruction lies in identifying objects, not in having a single one: there are strict test categories with a single object, such as the monoid of all self-maps of an infinite set and a countable monoid of primitive recursive functions.
\end{abstract}

\keywords{Monoid, Semigroup, homotopy types, classifying spaces, simplicial sets, test category
}
\subjclass[2020]{Primary 20M30, 55U35; Secondary 18F20, 18N40, 18B40, 20M05
}

\maketitle

\emergencystretch=1.5em


\section*{Introduction}

The idea that algebraic topology is the use of algebraic tools to study topology holds a truth, but only half of one. It is true that algebraic topology attaches groups, rings, modules\dots\ to spaces in a natural way; this is, in fact, what led to the creation of category theory \cite{kromer2007tool}. Moreover, these assignments have made it possible to answer important questions: to obtain obstructions to the existence of maps, to count fixed points or to study motion planning, to name a few. But all these approaches regard topology as the object of study and algebra as the instrument. To paraphrase a famous sentence, they ask what algebra can do for topology, but not what topology can do for algebra.

The history of topology at the service of algebra is shorter than the converse one, but it has some remarkable episodes. Perhaps the most prominent is the theorem of Nielsen and Schreier stating that every subgroup of a free group is free \cite{Schreier}. Its most elegant proof goes through topology. A free group $F$ is the fundamental group of a topological space, its classifying space $\B F$, and this space can be built as a $1$-dimensional complex, that is, as a graph. A subgroup of $F$ corresponds to a covering space, which is again a graph and has that subgroup as its fundamental group; being the fundamental group of a graph, the subgroup is free.

The important step in this argument is the functorial assignment of a topological space to each group, that is, the classifying space functor $\B\colon \mathrm{Group}\to\Top$, $G\mapsto\B G$, which goes precisely from algebraic objects to topological spaces. Once one realizes that such a bridge exists, it is natural to ask how much of topology lies at its far end. The first strong answer is the theorem of Kan and Thurston \cite{KT}: every connected space has the homology, with arbitrary local coefficients, of the classifying space of a group. At the level of homology, every connected space is essentially a group. But what about homotopy? It is almost entirely left outside: the space $\B G$ attached to a group has trivial homotopy groups from $\pi_2$ on.

The question then becomes how to attach topological spaces to other algebraic objects, so that they capture homotopy and not only homology. The answer came with the work of McDuff \cite{McDuff}, who proved that every path-connected space has the weak homotopy type of $\B M$ for some discrete monoid $M$. Her monoid has one generator for each simplex outside a maximal tree of a simplicial complex, invertible for the edges and \emph{idempotent} for the higher simplices, and the proof combines the methods of Kan--Thurston and of Baumslag--Dyer--Heller \cite{KT,BDH}. Moreover, $M$ can be chosen so that the inclusion of its maximal subgroup induces a homology equivalence: the group already sees the homology of Kan and Thurston, and what the non-invertible idempotents add is precisely the homotopy. This result has also been reformulated and extended by Fiedorowicz and Hurwitz \cite{Fied,Hurwitz}.

In this paper we follow this path, but with a kind of exchange of quantifiers. In McDuff's construction one is given a space and looks for a monoid that represents it; here we look for a monoid that represents all of them. Of course this makes no sense through the classifying space, which would give us a single space. What we shall do instead is to regard every topological space, at the level of its homotopy type, as a monoid acting on a set. The important point is that there is a monoid whose actions represent every homotopy type. In this paper we construct several monoids with this property, and we call them \emph{master monoids}.

In a way, this can be seen as a continuation of the spirit of Klein's Erlangen programme. There, spaces are studied through the transformations that preserve their structure, through the group that preserves their properties. This idea also motivates the passage from groups to categories, and in particular to the one-object categories that monoids are, since it is what led Eilenberg and Mac~Lane to introduce category theory \cite{EM45}. Spaces are understood through the way they are transformed into one another, and, for instance, all spaces are understood through their algebraic shadows under functors such as homology, cohomology or homotopy. Our proposal is an intermediate point: we study topological spaces as sets on which a single monoid acts.

At first sight this idea may seem strange, but the mystery disappears as soon as the trick is revealed. How can a single monoid represent every homotopy type through its actions? In the same way as a space is seen through its simplices, or through some similar combinatorial structure, that is, as a simplicial set, where the information lies in how the different cells are glued to obtain the space. The homotopy type of a space comes in this case as a functor to sets from some category. If we manage to \emph{collapse} that category into a single object, obtaining a monoid, and to do so while preserving the shape, we have the result we want. After all, a functor to sets from a monoid, seen as a category with a single object, is nothing but a set together with an action of the monoid on it.

The natural setting for this question is Grothendieck's theory of test categories, sketched in \emph{Pursuing Stacks} \cite{PS} and developed by Maltsiniotis and Cisinski \cite{Malt,Cis}. Grothendieck asked which small categories $\A$ of shapes can serve as a combinatorial foundation for homotopy theory, in the sense that the presheaves on $\A$, localized at the morphisms that induce homotopy equivalences between their categories of elements, recover the homotopy category of spaces. The simplicial category $\Delta$ is the founding example, but it is only one of many test categories, among them Joyal's category $\Theta$ \cite{CM}. From this point of view, our question asks whether the shapes of a test category can be compressed into a single object, and what survives of the test property when they are. Since a monoid is a category with a single object, it also leads to a natural question: whether there are test categories with only one object, that is, whether homotopy theory can be modelled by the actions of a single monoid on sets.

This journey is not without its tolls, and a substantial part of this paper is devoted to identifying them. Passing to the universal monoid preserves the ability to realize every homotopy type objectwise, but it generally loses the stronger ability to model the entire homotopy category. The main destination of the journey, together with its principal tolls, can be summarized as follows.

\begin{mainthm}
Let $M_\Delta=\Mon(\Dinj^{\op})$ be the universal monoid of the semisimplicial category, and $M_n=\Mon((\Dinj^{\leq n})^{\op})$ that of its truncations.
\begin{enumerate}
\item $M_\Delta$ is a master monoid: every CW homotopy type is of the form $S/\!/M_\Delta$ for some $M_\Delta$-set $S$.
\item For $n\geq1$, $M_n$ is finitely presented, with $n(n+3)/2$ generators, and the homotopy types of the form $S/\!/M_n$ are exactly those of dimension at most $n$.
\item If $\A$ is a weak test category, $\Mon(\A^{\op})$ is a master monoid. If moreover $\A$ has at least two objects, $\Mon(\A^{\op})$ is not a weak test category.
\item If $X$ is an infinite set, the monoid $\mathrm{End}(X)$ of all self-maps of $X$ is a strict test category; so is the countable monoid of primitive recursive functions $\mathbb N\to\mathbb N$.
\end{enumerate}
\end{mainthm}

Statement (4) shows that the toll in (3) is the price of identifying objects, not the price of having a single object (Proposition~\ref{prop:test-monoid} and Corollary~\ref{cor:test-monoid-pr}). The test structure of $\mathrm{End}(X)$ does not come from identifying the objects of a test category, but from keeping them alive as idempotents.

The paper is organized as follows. In Section~\ref{sec:categories-topological-combinatorial} we recall how categories act on sets and how they determine homotopy types. In Section~\ref{sec:universal-monoid} we introduce the universal monoid of a category, its normal forms and the free action, which we read as an \emph{Aufhebung}. Section~\ref{sec:coref} contains the coreflection theorem, which is the heart of the paper. In Section~\ref{sec:semisimplicial-master-monoid} we study the semisimplicial master monoid and its truncations; the exactness of the truncated monoids is the only place where we need to control arbitrary, non-free actions. In Section~\ref{sec:test-categories-master-monoids} we compare master monoids with Grothendieck's test categories, show that a single monoid can nevertheless model homotopy theory when shapes are split instead of identified, and end with a reflection on what has been gained and what has been lost.

\subsection*{Use of artificial intelligence} 
In preparing this paper the author used Claude Opus~5.5 (Anthropic) to improve the writing, to translate the text from Spanish into English, to review the logic of the arguments, and to check whether relevant literature on the subject had been omitted. The diagrams of Figure~\ref{fig:two-edges}, Example~\ref{ex:sphere}, Figure~\ref{fig:kron} and the figure in Remark~\ref{rem:point} are due to Claude. The author takes full responsibility for the content of the paper.

\section{From categories and actions to homotopy types}
\label{sec:categories-topological-combinatorial}

The categories that usually appear in mathematics are large: their collections of objects and morphisms are proper classes. In this paper we adopt the opposite point of view. All categories will be small, and they will be regarded as mathematical objects in their own right; we shall study their combinatorial structure and the homotopy type they determine.

A particularly useful way of describing a small category is by generators and relations, in analogy with presentations of semigroups, monoids and groups. Following \cite{Poly}, a categorical presentation consists of a quiver $Q=(Q_0,Q_1,s,t)$ and a family $R$ of relations between parallel paths; the presented category is the free category on $Q$ divided by the congruence generated by $R$. We shall be especially interested in categories with a single object, that is, in monoids. A monoid $M$ is a category whose only object we denote by $*$ and whose endomorphisms are the elements of $M$, with composition given by the product. With our convention, $gf=g\circ f$ means that $f$ acts first and then $g$. Groups are precisely the one-object categories in which every morphism is invertible, and a presentation of a monoid as a category, with a single vertex and one loop for each generator, is the same as an ordinary presentation of the monoid.

\subsection{Categories acting on sets}
\label{subsec:categories-acting-sets}

An action of a monoid or a group on a set can be read as a functor from the corresponding one-object category to $\Set$. This suggests studying, for a small category $\CC$, functors $X\colon\CC\to\Set$. We may think of a category as a monoid with several objects, and of such a functor as an action on a set with several types.

We work with covariant functors. A \emph{$\CC$-set} is a functor $X\colon\CC\to\Set$. An element $x\in X(c)$ has \emph{type} $c$, and each morphism $f\colon c\to c'$ transforms elements of type $c$ into elements of type $c'$ through the map $X(f)\colon X(c)\to X(c')$. The Grothendieck construction attaches to $X$ its \emph{category of elements} $\int_{\CC}X$, whose objects are the pairs $(c,x)$ with $x\in X(c)$, a morphism $(c,x)\to(c',x')$ being a morphism $f\colon c\to c'$ of $\CC$ with $X(f)(x)=x'$. When there is no risk of confusion we simply write $x$ for $(c,x)$. The projection $\int_{\CC}X\to\CC$ is a discrete opfibration, and in fact the Grothendieck construction gives an equivalence between $\CC$-sets and discrete opfibrations over $\CC$. This reading is also the basis of the functorial model of databases: a category is a schema, a functor to $\Set$ is an instance, and functors between schemas induce data migration operations \cite{SevenSketches,SpivakDataMigration,SpivakDatabaseQueries}.

\subsubsection{Monoids acting on sets}
\label{subsubsec:monoid-actions}

Let $M$ be a monoid. An \emph{$M$-set} is a functor $S\colon M\to\Set$, that is, a set $S$ with a left action $(m,s)\mapsto m\cdot s$. Its category of elements is the \emph{translation category} $M\ltimes S$: its objects are the elements of $S$, and
\[
\Hom_{M\ltimes S}(s,t)=\{m\in M\mid m\cdot s=t\},
\]
with composition induced by the product of $M$. In particular, out of every object $s$ comes exactly one morphism labelled by each element $m\in M$, with target $m\cdot s$. We call the classifying space $S/\!/M=B(M\ltimes S)$ the \emph{homotopy quotient} of the action. When $M$ is a group, $M\ltimes S$ is a groupoid, the \emph{action groupoid} of the action, and $S/\!/M$ is the usual homotopy quotient of a set by a group (Example~\ref{ex:groups}). For the one-point set we get $*/\!/M=BM$, and for the left regular action $M/\!/M\simeq*$, because the unit $1\in M$ is an initial object of $M\ltimes M$.

\begin{example}[Groups]
\label{ex:groups}
If a group $G$ acts on a set $S$, then $G\ltimes S$ is a groupoid whose connected components are the orbits, and the automorphism group of an object $s$ is its stabilizer $G_s$. Hence $S/\!/G\simeq\coprod_{[s]\in S/G}BG_s$. In particular, homotopy quotients of group actions on sets are disjoint unions of aspherical spaces.
\end{example}

\begin{example}[One transformation]
\label{ex:N}
An $\mathbb N$-set is a set $S$ with an endomorphism $\varphi\colon S\to S$. The category $\mathbb N\ltimes S$ is the free category on the functional graph $s\to\varphi(s)$, and the classifying space of a free category has the homotopy type of its generating graph. Each component of a functional graph contains at most one cycle, so each component of $S/\!/\mathbb N$ is contractible or a circle:
\[
\begin{array}{c@{\hspace{2.2em}}c@{\hspace{2.2em}}c}
\begin{tikzcd}[column sep=small]
0 \arrow[r] & 1 \arrow[r] & 2 \arrow[ll,bend left=40]
\end{tikzcd}
&
\begin{tikzcd}[column sep=small]
0 \arrow[r] & 1 \arrow[r] & 2 \arrow[r] & \cdots
\end{tikzcd}
&
\begin{tikzcd}[column sep=small]
x \arrow[r] & y \arrow[r,bend left=30] & z \arrow[l,bend left=30]
\end{tikzcd}
\\[4pt]
S=\mathbb Z/3,\ \simeq S^1 & S=\mathbb N,\ \simeq * & S=\{x,y,z\},\ \simeq S^1
\end{array}
\]
\end{example}

\begin{example}[Two transformations]
\label{ex:F2}
Let $F_2=\langle s,t\rangle$ be the free monoid on two generators. For every $F_2$-set $S$, the category $F_2\ltimes S$ is free on the graph with two outgoing arrows at each vertex, $x\to s\cdot x$ and $x\to t\cdot x$. Hence $S/\!/F_2$ has the homotopy type of a graph. If a finite component has $k$ elements, its graph has $k$ vertices and $2k$ edges, and it is a bouquet of $2k-k+1=k+1$ circles.
\end{example}

\subsection{Categories and topological spaces}
\label{subsec:categories-topological-spaces}

Every small category determines functorially a topological space, its classifying space. This extends the classical construction $G\mapsto BG$ for groups and allows us to study categories by homotopical methods.

\subsubsection{The nerve and the classifying space}

The \emph{nerve} of a small category $\CC$ is the simplicial set $N\CC\colon\Delta^{\op}\to\Set$ whose $n$-simplices are the strings of $n$ composable morphisms $x_0\to x_1\to\cdots\to x_n$; equivalently, $N_n\CC=\Fun([n],\CC)$, where $[n]=\{0<1<\cdots<n\}$ is regarded as a category. The face operators compose two consecutive arrows or drop one of the ends of the string, and the degeneracies insert identities.

\begin{definition}
The \emph{classifying space} of a small category $\CC$ is the geometric realization of its nerve, $B\CC=|N\CC|$.
\end{definition}

The construction $B\colon\Cat\to\Top$ is functorial, and a natural transformation between two functors induces a homotopy between the corresponding maps. In particular, an adjunction $F\colon\CC\rightleftarrows\DD\colon G$ induces homotopy equivalences $B\CC\simeq B\DD$, and a category with an initial or a terminal object has a contractible classifying space \cite{LEE,LEE2}. The construction preserves products, $B(\EE\times\FF)\cong B\EE\times B\FF$, but not ordinary colimits in general. The right relation between categorical constructions and homotopy colimits is given by Thomason's theorem.

\subsubsection{The Grothendieck construction and Thomason's theorem}

The category of elements of a functor to sets is the discrete case of a more general construction. Let $F\colon\CC\to\Cat$ be a diagram of small categories. Its \emph{Grothendieck construction} $\int_{\CC}F$ has as objects the pairs $(c,x)$ with $c\in\Ob\CC$ and $x\in\Ob F(c)$; a morphism $(c,x)\to(c',x')$ is a pair $(f,u)$ formed by a morphism $f\colon c\to c'$ of $\CC$ and a morphism $u\colon F(f)(x)\to x'$ of $F(c')$. The following result is the fundamental tool that allows us to read the Grothendieck construction as a categorical procedure of homotopical assembly.

\begin{theorem}[Thomason]
\label{thm:thomason}
Let $\CC$ be a small category and $F\colon\CC\to\Cat$ a functor. There is a natural weak equivalence of simplicial sets $\hocolim_{c\in\CC}NF(c)\to N\!\left(\int_{\CC}F\right)$, and hence a natural homotopy equivalence
\[
\hocolim_{c\in\CC}BF(c)\simeq B\!\int_{\CC}F.
\]
\end{theorem}

This is Thomason's homotopy colimit theorem \cite[Theorem~1.2]{Th}. It says that the Grothendieck construction is a strict categorical model of the homotopy colimit of a diagram of categories. Unlike the ordinary colimit, which may lose the information about how the objects of the diagram are assembled, the Grothendieck construction keeps both the categories $F(c)$ and the morphisms that describe how they vary along $\CC$. Regarding each $X(c)$ as a discrete category, for a $\CC$-set we obtain
\[
B\!\int_{\CC}X\simeq\hocolim_{c\in\CC}X(c),
\]
where the right-hand side is the homotopy colimit of the diagram of discrete spaces underlying $X$.

This must not be confused with the ordinary colimit $\colim_{\CC}X$, which identifies elements along the morphisms of $\CC$ but forgets the homotopical information contained in the different ways of making those identifications. For the constant functor with value a point we have $\int_{\CC}*\cong\CC$, and hence $\hocolim_{\CC}*\simeq B\CC$, whereas the ordinary colimit of a nonempty connected $\CC$ is just a point. Even a diagram made only of points can produce a space with nontrivial homotopy. Thomason's theorem is also what will allow us to read the coreflection theorem homotopically, as a statement about Kan extensions (Corollary~\ref{cor:kan}).

\begin{definition}
\label{def:realization-spectrum}
The \emph{realization spectrum} of a small category $\CC$ is the class
\[
\T(\CC)=\left\{\left[B\!\int_{\CC}X\right]\ \middle|\ X\colon\CC\to\Set\right\}
\]
of homotopy types represented by the categories of elements of $\CC$-sets; by Thomason's theorem, these are the homotopy types $\hocolim_{\CC}X$, with the values of $X$ regarded as discrete spaces.
\end{definition}

The spectrum $\T(\CC)$ measures how much topology can be encoded by actions of $\CC$ on sets with types. It is not only the homotopy type of $B\CC$, which corresponds to the constant functor with value a point, but the whole collection of homotopy types obtained as the diagram $X$ varies.

\begin{example}
\label{ex:sorts}
\begin{enumerate}[leftmargin=*]
\item \emph{Graphs.} Let $\mathbf K=\{E\rightrightarrows V\}$ be the Kronecker category, with arrows $s,t\colon E\to V$. A functor $X\colon\mathbf K\to\Set$ is a graph, or quiver, with $X(E)$ the set of edges, $X(V)$ the set of vertices, and source and target maps $X(s),X(t)\colon X(E)\to X(V)$. The category $\int_{\mathbf K}X$ is the incidence category of the graph, its nerve gives the barycentric subdivision, and $B\!\int_{\mathbf K}X\simeq|X|$. For the graph with one vertex $v$ and one loop $e$, the category of elements is $e\rightrightarrows v$, and $B\!\int_{\mathbf K}X\simeq S^1$.
\item \emph{Semisimplicial sets.} A functor $X\colon\Dinj^{\op}\to\Set$ is a semisimplicial set. Its category of elements is its category of simplices, and there is a natural homotopy equivalence $B\!\int_{\Dinj^{\op}}X\simeq|X|$, which can also be seen as the discrete case of Thomason's theorem.
\item \emph{One-object categories.} If $G$ is a group, Example~\ref{ex:groups} shows that $\T(G)$ consists of the disjoint unions of classifying spaces of subgroups of $G$. Example~\ref{ex:N} shows that $\T(\mathbb N)$ consists of disjoint unions of points and circles, and Example~\ref{ex:F2} that $\T(F_2)$ is contained in the class of homotopy types of graphs. We shall see in Example~\ref{ex:spectra} that it contains all of them.
\end{enumerate}
\end{example}

The question raised in the introduction can now be stated precisely: is it possible to forget the types, that is, the objects of a category, and to replace a category with several objects by a monoid without changing its realization spectrum? The coreflection theorem will answer affirmatively, in the sense that $\T(\CC)\subseteq\T(\Mon(\CC))$. More concretely, it will show how to turn each $X\colon\CC\to\Set$ into a $\Mon(\CC)$-set whose homotopy quotient has the homotopy type of $B\int_{\CC}X$.

\section{The universal monoid of a category}
\label{sec:universal-monoid}

Given a category, let us see how to attach a monoid to it in a \emph{natural} way. Recall that not only is every monoid a small category, but every morphism of monoids is exactly a functor. Hence we have a functor $\mathrm{Mon}\hookrightarrow\Cat$, and its left adjoint will give us the object we are looking for. Throughout the section we present the ideas through examples, and we close it with an example that gathers all of them (Section~\ref{subsec:example-two-edges}).

\begin{defi}
The \emph{universal monoid} $\Mon(\C)$ of a small category $\C$ is the pushout in $\Cat$
\[
\begin{tikzcd}
\Ob\C \arrow[r] \arrow[d] & \C \arrow[d,"\iota"]\\
{*} \arrow[r] & \Mon(\C),
\end{tikzcd}
\]
where $\Ob\C$ denotes the discrete category on the objects of $\C$. Equivalently, $\Mon$ is left adjoint to the inclusion of monoids into categories, with unit $\iota\colon\C\to\Mon(\C)$. Hence a functor from $\C$ to a monoid $N$ is the same as a homomorphism $\Mon(\C)\to N$.
\end{defi}

Using the description of colimits of presented categories in \cite[\S3.1]{Poly}, we see that $\Mon(\CC)$ is the monoid generated by all the morphisms of $\CC$ subject to $[g][f]=[g\circ f]$ for every composable pair, where identities become the unit. More explicitly, if a small category is presented by a quiver and relations, $\C=\langle Q\mid R\rangle$, then $\Mon(\C)=\langle Q_1\mid R\rangle$, where $Q_1$ is the set of arrows.

\begin{example}
\label{ex:universal-monoids}
\begin{enumerate}[leftmargin=*]
\item \emph{Parallel arrows.} In the Kronecker category $\mathbf K=\{E\rightrightarrows V\}$, with arrows $s,t\colon E\to V$, no two non-identity arrows are composable, so no relation is imposed and $\Mon(\mathbf K)=\langle s, t \rangle=F_2$ is the free monoid on two generators. In the same way $\Mon(\{0\to1\})=\mathbb N$, and the universal monoid of the free category on a quiver is the free monoid on its arrows.
\item \emph{A category with relations.} Let $\mathbf P$ be the category with objects $x,y,z$ generated by two arrows $f_1,f_2\colon x\to y$ and two arrows $g_1,g_2\colon y\to z$, subject to
\[
\begin{tikzcd}[column sep=large]
x \arrow[r, shift left, "f_1"] \arrow[r, shift right, "f_2"'] & y \arrow[r, shift left, "g_1"] \arrow[r, shift right, "g_2"'] & z
\end{tikzcd}
\qquad g_1f_1=g_2f_2,\qquad g_2f_1=g_1f_2 .
\]
It has six non-identity arrows: $f_1,f_2,g_1,g_2$ and the two composites $h=g_1f_1=g_2f_2$ and $k=g_2f_1=g_1f_2$. Its universal monoid is
\[
\Mon(\mathbf P)=\langle f_1,f_2,g_1,g_2\mid g_1f_1=g_2f_2,\ g_2f_1=g_1f_2\rangle .
\]
\end{enumerate}
\end{example}

\subsection{Normal forms}

When we study a mathematical object we like to have a manageable way of working with it. For the monoid attached to a category we can use a presentation of the category to obtain a presentation of the monoid; but for a monoid a presentation is only a first step, and we would also like to have normal forms.

An element of $\Mon(\CC)$ is a product $u_1\cdots u_k$ of non-identity arrows of $\CC$, which we call \emph{letters}; the letter $u_k$ acts first. A word is \emph{reduced} when no two consecutive letters are composable, that is, when $\dom u_i\neq\cod u_{i+1}$ for all $i$. Let us see that this gives the normal forms.

\begin{lemma}[Normal forms]\label{lem:nf}
For every small category $\C$, each element of $\Mon(\C)$ is represented by a unique reduced word. In particular, $\iota$ is faithful and injective on non-identity arrows, and $\iota^{-1}(1)$ consists of the identities.
\end{lemma}

\begin{proof}
Consider the $2$-polygraph with one $0$-cell, one $1$-cell for each non-identity arrow of $\C$, and one $2$-cell $uv\Rightarrow u\circ v$ for each composable pair, where the target is the empty word when $u\circ v$ is an identity. It presents $\Mon(\C)$, and it terminates because every rule decreases length. Its critical branchings are the words $uvw$ with $(u,v)$ and $(v,w)$ composable. Both reductions end in $u\circ v\circ w$ by associativity; when one of the intermediate composites is an identity, say $u\circ v=\mathrm{id}$, one side gives $w$ and the other gives $u\circ(v\circ w)=w$, and similarly in the other cases. By the critical branching lemma \cite[Lemma~4.3.7]{Poly} and Newman's lemma \cite[Lemma~1.3.21]{Poly}, the polygraph is convergent, and its normal forms are the reduced words.
\end{proof}

In view of this lemma we shall write simply $f$, instead of $\iota(f)$, for the element of $\Mon(\CC)$ represented by a morphism $f$ of $\CC$, the identities becoming the unit $1$. This is a mild abuse of notation: a non-identity morphism $f$ is at the same time an arrow of $\CC$ and a letter of the universal monoid, and the context will always make clear which role it plays. For the category $\mathbf P$ of Example~\ref{ex:universal-monoids}, the letters are the six arrows $f_1,f_2,g_1,g_2,h,k$, and a word is reduced exactly when no letter $g_i$ is immediately followed by a letter $f_j$; for instance, $g_2f_1$ is not reduced, and its normal form is the one-letter word $k$. A presentation needs only four generators, but the normal forms need all six arrows.

Let us now see a lemma that will be the key to the main theorem. In a category of elements $\int X$ over $\CC$, the arrows going out of an object $(c,x)$ are exactly those going out of $c$ in $\CC$. For a monoid $M$, on the contrary, out of every object of any translation category come as many arrows as there are elements of $M$. Let us see what this means for $\Mon(\CC)$ when we decompose the monoid with respect to a fixed object $c$ of $\CC$. For an object $c$, let $W_c$ be the set of reduced words that are either empty or whose first-acting letter $u_k$ satisfies $\dom u_k\neq c$. These are the words that cannot absorb anything sitting at $c$.

\begin{lemma}[Unique factorization]\label{lem:fact}
For every object $c$ of $\C$, the map
\[
\coprod_dW_d\times\C(c,d)\longrightarrow\Mon(\C),\qquad (w,g)\longmapsto wg,
\]
is a bijection.
\end{lemma}

\begin{proof}
Let $m$ have reduced word $u_1\cdots u_k$. If $k\ge1$ and $\dom u_k=c$, then $m=(u_1\cdots u_{k-1})\,u_k$, and $u_1\cdots u_{k-1}\in W_{\cod u_k}$ because the word is reduced. Otherwise $m\in W_c$, and $m$ is the image of $(m,\mathrm{id}_c)$. Conversely, if $m=wg$ with $g\neq\mathrm{id}$ and $w\in W_{\cod g}$, then the word $wg$ is reduced, so it is the normal form of $m$ and $g$ is its first-acting letter; if $g=\mathrm{id}$, then $m=w\in W_c$. The two cases exclude each other.
\end{proof}

In other words, think of an element $x$ of type $c$. Everything in $\Mon(\CC)$ that may act on it is a word whose first-acting letter either goes out of $c$, and can then be absorbed by $x$, or does not, and then the word belongs to $W_c$; this last case is the one in which the factor in $\CC(c,c)$ is the identity. Let us now see how to build sets on which $\Mon(\CC)$ acts.

\subsection{\emph{Aufhebung}: restriction and free actions}
\label{subsec:restriction-free-actions}

The relation between a category $\CC$ and its universal monoid is not just a construction, not even a functorial one. It comes with a canonical functor $\iota\colon\CC\to\Mon(\CC)$, which sends all objects of $\CC$ to the only object of $\Mon(\CC)$ and each morphism $f$ to the element $f$ it represents in the universal monoid. This functor will serve as a bridge from a functor $X\colon\CC\to\Set$ to a set on which $\Mon(\CC)$ acts, and it will allow us to replace the category of elements $\int_{\CC}X$ by the translation category of an action of the universal monoid without changing its homotopy type. When we forget the different types given by the objects of $\CC$ and replace them by a single total action, new elements appear; but, as we shall see, they are organized in branches that do not alter the original homotopy type.

Recall first the general categorical situation. Let $F\colon\CC\to\DD$ be a functor between small categories and $\A$ a complete and cocomplete category. Precomposition with $F$ defines the restriction functor $F^*\colon\A^{\DD}\to\A^{\CC}$, $F^*(Y)=Y\circ F$, which has a left and a right adjoint, the left and right Kan extensions $F_!=\Lan_F$ and $F_*=\Ran_F$:
\[
\begin{tikzcd}[column sep=6em]
\A^{\CC}
\ar[r, bend left=35, "F_!"]
\ar[r, bend right=35, swap, "F_*"]
&
\A^{\DD}
\ar[l, "F^*" description]
\end{tikzcd}
\qquad
F_!\dashv F^*\dashv F_*.
\]
Thus the restriction of a diagram has two associated universal constructions, a free extension given by $F_!$ and a cofree extension given by $F_*$; for the general theory of adjunctions and Kan extensions see \cite[Chapters~4 and~6]{RiehlCategoryTheoryContext}. A triple of adjoints of this form also appears in the study of the levels of a presheaf topos: a full inclusion of small categories induces a restriction functor between the corresponding presheaf categories, together with its left and right Kan extensions \cite[Section~1]{KennettRiehlRoyZaks2011}. In that context this configuration is related to Lawvere's \emph{unity and identity of opposites}, and the term \emph{Aufhebung}, introduced by Lawvere in allusion to Hegel, expresses the way in which a structure can be at the same time abolished and preserved \cite{LawvereFutureCategoryTheory}.

Let us apply this to the canonical functor $\iota\colon\CC\to M$, with $M=\Mon(\CC)$, and to $\A=\Set$. We obtain the triple of adjoints
\[
\begin{tikzcd}[column sep=6em]
\Set^{\CC}
\ar[r, bend left=35, "\iota_!"]
\ar[r, bend right=35, swap, "\iota_*"]
&
\Set^M
\ar[l, "\iota^*" description]
\end{tikzcd}
\qquad
\iota_!\dashv\iota^*\dashv\iota_*.
\]
The restriction $\iota^*$ sends an $M$-set $S$ to the $\CC$-set which has the same set $S$ at every object of $\CC$, an arrow $f\colon c\to c'$ acting through the letter $f\in M$. After restricting the action, each element of $S$ belongs simultaneously to all the types given by the objects of $\CC$. Its left adjoint, the left Kan extension, is the functor we shall mainly be interested in. For $X\colon\CC\to\Set$ it is given by
\[
\iota_!X=M\otimes_{\CC}X=\int^{c\in\CC}M(\iota c,*)\times X(c),
\]
the universal $M$-set equipped with a natural transformation $X\to\iota^*\iota_!X$. Let us describe it explicitly.

\begin{proposition}[The free action]
\label{prop:free}
Let $M=\Mon(\CC)$ and $X\colon\CC\to\Set$. Every element of $\iota_!X$ can be written uniquely as $[w,x]=w\cdot[1,x]$, with $x\in X(c)$ and $w\in W_c$, so that
\[
\iota_!X\cong\coprod_{c\in\Ob\CC}W_c\times X(c).
\]
Moreover, if $m\in M$ factors as $m=w'g$ with $g\colon c\to c'$ and $w'\in W_{c'}$, then $m\cdot[1,x]=[w',X(g)x]$.
\end{proposition}

\begin{proof}
By the co-Yoneda lemma and Lemma~\ref{lem:fact},
\[
\begin{aligned}
M\otimes_{\CC}X&=\int^{c\in\CC}M(\iota c,*)\times X(c)\cong\int^{c\in\CC}\Bigl(\coprod_{d\in\Ob\CC}W_d\times\CC(c,d)\Bigr)\times X(c)\\
&\cong\coprod_{d\in\Ob\CC}W_d\times\int^{c\in\CC}\CC(c,d)\times X(c)\cong\coprod_{d\in\Ob\CC}W_d\times X(d),
\end{aligned}
\]
and if $m=w'g$ is the unique factorization of $m$ with respect to $c$, the class of $(m,x)$ corresponds to $(w',X(g)x)$. This gives both the description and the uniqueness of the normal form $[w,x]$.
\end{proof}

Thus $\iota_!X$ is the $M$-set presented by the elements of $X$ as generators and the relations $f\cdot[1,x]=[1,X(f)x]$, for $f\colon c\to c'$ and $x\in X(c)$. When a letter acts on a \emph{root} $[1,x]$, it is absorbed into the structure of $X$ whenever it represents an arrow whose domain is the type of $x$. If it cannot be absorbed, it starts a new \emph{branch}.

This description explains our use of the word \emph{Aufhebung}. In the theory of levels of a topos, this terminology expresses a relation in which a structure is overcome and preserved at the same time \cite[Section~1]{KennettRiehlRoyZaks2011}. In our case, the types determined by the objects of $\CC$ are abolished as explicit classes when we pass to the universal monoid, but the information they carried does not disappear completely: it remains encoded in the operations that each element can absorb. The passage from $X$ to $\iota_!X$ forgets the explicit classification of elements by types, but keeps it in an operational form.

\begin{remark}[The cofree action]
\label{rem:cofree}
The restriction also has a right adjoint, $\iota_*X=\Hom_{\Set^{\CC}}(\iota^*M,X)$, where $M$ is the regular $M$-set and the action of $M$ is induced by its right action on itself. The elements of $\iota_*X$ can be read as the possible \emph{behaviours} of an element under all the words of $M$. While the left extension adds freely the results of the operations that cannot be performed inside $X$, the right extension prescribes simultaneously compatible answers for all of them. Unlike $\iota_!$, it does not preserve the homotopy type of the category of elements. For the arrow category $\CC=\{0\to1\}$ and the terminal functor $X=(*\to*)$, the category $\int_{\CC}X$ has a terminal object and $B\!\int_{\CC}X\simeq*$; but $\Mon(\CC)\cong\mathbb N$, $\iota_*X$ is a point with the trivial action, and $B(\mathbb N\ltimes\iota_*X)\cong B\mathbb N\simeq S^1$. Homotopy invariance is therefore a property of the free side of the triple, of $\iota_!$, and not of the cofree extension $\iota_*$.
\end{remark}

\subsection{An example}
\label{subsec:example-two-edges}

Let us gather everything seen so far in a single example, on the Kronecker category $\mathbf K=\{E\rightrightarrows V\}$, whose universal monoid is the free monoid $F_2=\langle s,t\rangle$ (Example~\ref{ex:universal-monoids}). Recall that a $\mathbf K$-set $X$ is a quiver, with set of edges $X(E)$, set of vertices $X(V)$, source map $X(s)$ and target map $X(t)$.

\begin{figure}[!htbp]
\[
\begin{array}{c@{\hspace{3em}}c}
\begin{array}{c}
\begin{tikzcd}[column sep=large]
a \arrow[r, shift left=1.5, "e_1"] \arrow[r, shift right=1.5, "e_2"'] & b
\end{tikzcd}
\\[16pt]
\begin{tikzcd}[row sep=small, column sep=small]
 & & a & & \\
e_1 \arrow[urr,"s"] \arrow[drr,"t"'] & & & & e_2 \arrow[ull,"s"'] \arrow[dll,"t"] \\
 & & b & &
\end{tikzcd}
\end{array}
&
\begin{tikzcd}[row sep=small, column sep=small]
 & {[s,a]} & & {[t,a]} & \\
 & & a \arrow[ul,dashed,"s"'] \arrow[ur,dashed,"t"] & & \\
e_1 \arrow[urr,"s"] \arrow[drr,"t"'] & & & & e_2 \arrow[ull,"s"'] \arrow[dll,"t"] \\
 & & b \arrow[dl,dashed,"s"] \arrow[dr,dashed,"t"'] & & \\
 & {[s,b]} & & {[t,b]} &
\end{tikzcd}
\\[6pt]
\text{(a) } X \text{ and } \int_{\mathbf K}X & \text{(b) } F_2\ltimes\iota_!X
\\[14pt]
\begin{tikzcd}[column sep=small]
 & 2 \arrow[dr,"e_2"] & \\
1 \arrow[ur,"e_1"] & & 3 \arrow[ll,"e_3"]
\end{tikzcd}
&
\begin{tikzcd}[row sep=large, column sep=normal]
 & e_1 \arrow[dl,"s"'] \arrow[dr,"t"] & \\
1 & & 2 \\
e_3 \arrow[u,"t"] \arrow[dr,"s"'] & & e_2 \arrow[u,"s"'] \arrow[dl,"t"] \\
 & 3 &
\end{tikzcd}
\\[6pt]
\text{(c) the quiver } \iota^*S & \text{(d) } \int_{\mathbf K}\iota^*S
\end{array}
\]
\[
\begin{array}{c}
\begin{tikzcd}[column sep=small]
 & 2 \arrow[dr,"t"] \arrow[loop above,"s"] & \\
1 \arrow[ur,"t"] \arrow[loop left,"s"] & & 3 \arrow[ll,"t"] \arrow[loop right,"s"]
\end{tikzcd}
\\[6pt]
\text{(e) } F_2\ltimes S
\end{array}
\]
\caption{The example of Section~\ref{subsec:example-two-edges}. (a) The quiver $X$, with two vertices $a,b$ and two parallel edges $e_1,e_2$ from $a$ to $b$ (top), and its category of elements, a crown (bottom). (b) The translation category of the free action $\iota_!X$: the edges absorb both letters (solid arrows), while from each vertex grows a binary tree of branches (dashed arrows), drawn up to its first level. (c) The quiver $\iota^*S$ obtained by restricting the action of $F_2$ on $\{1,2,3\}$ in which $s$ is the identity and $t$ is the cyclic permutation $1\mapsto2\mapsto3\mapsto1$. (d) Its category of elements, a hexagonal crown. (e) The translation category of the action itself.}
\label{fig:two-edges}
\end{figure}
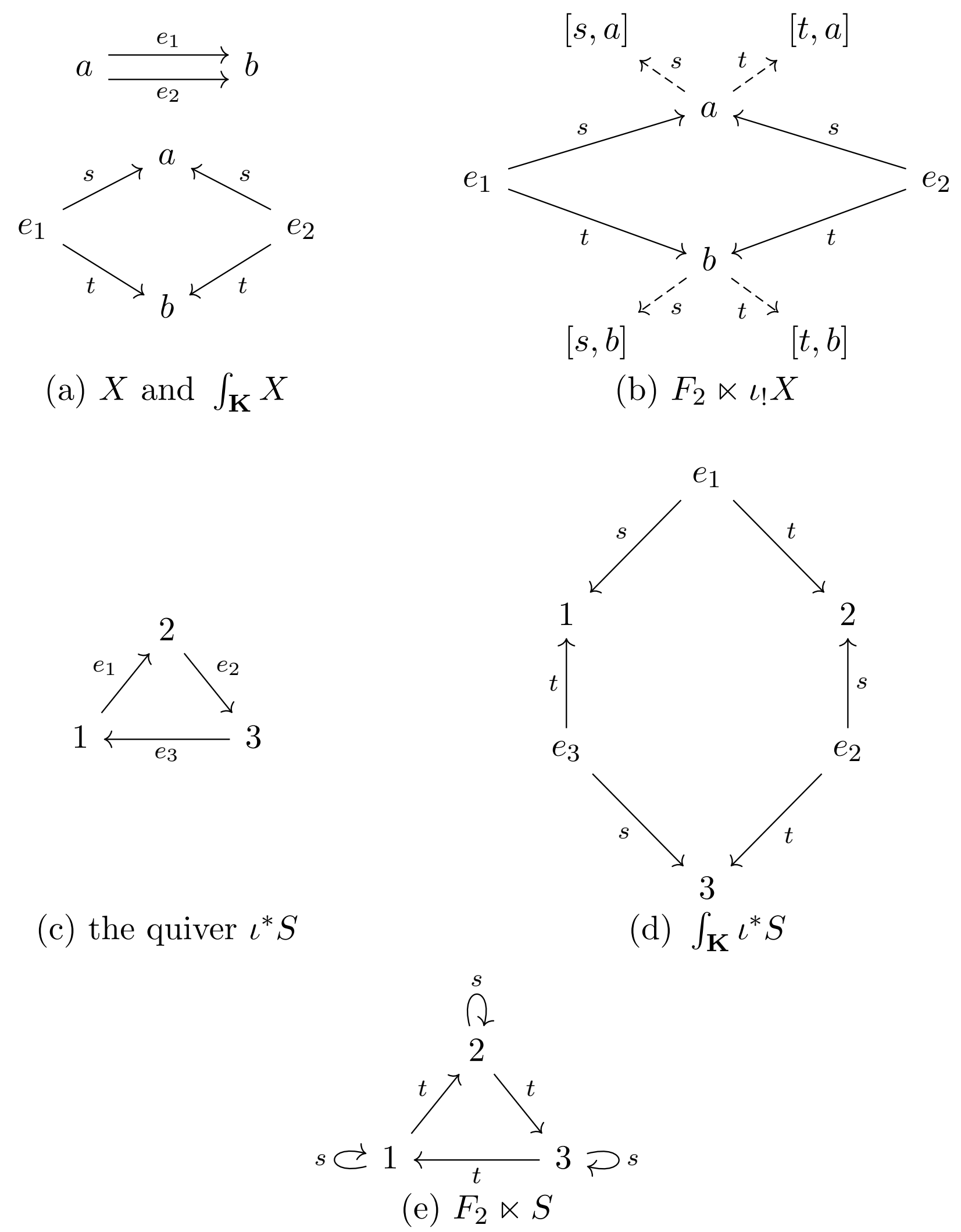

\begin{example}
\label{ex:two-edges}
Let $X$ be the quiver with two vertices $a,b$ and two edges $e_1,e_2$ from $a$ to $b$, that is, $X(E)=\{e_1,e_2\}$, $X(V)=\{a,b\}$, $X(s)(e_i)=a$ and $X(t)(e_i)=b$.

\emph{The category of elements.} The category $\int_{\mathbf K}X$ has the four objects $e_1,e_2,a,b$ and an arrow from each edge to each of its ends; it is the poset known as the crown, Figure~\ref{fig:two-edges}(a). Its classifying space is a square, the barycentric subdivision of the graph, so $B\!\int_{\mathbf K}X\simeq S^1$.

\emph{The free action.} Both letters $s,t$ have domain $E$. Hence $W_E=\{1\}$, since every nonempty reduced word acts first by a letter with domain $E$, and $W_V=F_2$. By Proposition~\ref{prop:free},
\[
\iota_!X\cong\{e_1,e_2\}\sqcup\big(F_2\times\{a,b\}\big),\qquad
s\cdot e_i=a,\quad t\cdot e_i=b,\quad u\cdot[w,v]=[uw,v],
\]
for $v\in\{a,b\}$ and $u,w\in F_2$. The edges absorb both letters, whereas the vertices absorb none, and from each of them grows a binary tree of branches, Figure~\ref{fig:two-edges}(b). The core is the crown, and the coreflection theorem of the next section will show that the branches do not change the homotopy type: $\iota_!X/\!/F_2\simeq S^1$.

\emph{Restricting an action.} Conversely, let $S=\{1,2,3\}$ with $s$ acting as the identity and $t$ as the cyclic permutation $1\mapsto2\mapsto3\mapsto1$, written $231$. The $\mathbf K$-set $\iota^*S$ has $S$ both as its set of edges and as its set of vertices, and the edge $e_i$ goes from $s\cdot i=i$ to $t\cdot i$. It is the oriented triangle of Figure~\ref{fig:two-edges}(c), and its category of elements is the hexagonal crown of Figure~\ref{fig:two-edges}(d), so $B\!\int_{\mathbf K}\iota^*S\simeq S^1$. The action itself sees more: $F_2\ltimes S$ is free on the graph of Figure~\ref{fig:two-edges}(e), with three vertices and six edges, and $S/\!/F_2$ is a bouquet of four circles. The three extra circles appear because each element of $S$ is, after restriction, both an edge and a vertex, and the action does not distinguish these two copies; Lemma~\ref{lem:arbitrary-actions} will make this precise.
\end{example}

\section{The coreflection theorem}\label{sec:coref}

In this section we prove the main result towards the construction of a master monoid. We shall see how, given a category $\CC$, the homotopy type of the category of elements $\int_{\CC}X$ is represented by the homotopy type of $\iota_!X$ as a $\Mon(\CC)$-set. From now on $\CC$ is a small category, $X\colon\CC\to\Set$ and $M=\Mon(\CC)$.

Let us first see what is at stake. The passage from $X$ to $\iota_!X$ is not innocent: every operation that an element cannot absorb opens a new branch, and these branches are usually infinite, even when $X$ is finite. In the category of elements, an object $(c,x)$ only emits the arrows of $\CC$ that go out of its type $c$; in the translation category $M\ltimes\iota_!X$, every object emits one arrow for each element of $M$. The second category is therefore, in general, enormously larger than the first, and there is no a priori reason why both should have the same homotopy type. The theorem below says that all this new material is homotopically invisible. The category of elements sits inside the translation category as a full subcategory, which we call the \emph{core}, and the branches retract onto their roots through an adjunction. Recall from Section~\ref{sec:categories-topological-combinatorial} that an adjunction induces a homotopy equivalence of classifying spaces; this is all the homotopy theory we need.

\begin{theorem}[Coreflection theorem]\label{thm:B}
Let $M=\Mon(\C)$ and $X\colon\C\to\Set$. The functor $V\colon\int X\to M\ltimes\iota_!X$, given by $V(c,x)=[1,x]$ and $V(f)=f$, is fully faithful and has a right adjoint $R$, given on objects by $R[w,x]=(c,x)$. The counit at $[w,x]$ is the arrow $w\colon[1,x]\to[w,x]$. Hence
\[
B\!\textstyle\int_\C X\ \simeq\ \iota_!X/\!/\Mon(\C),
\]
and $\T(\C)\subseteq\T(\Mon(\C))$.
\end{theorem}

\begin{proof}
Let $(c,x)\in\int X$ and $s=[w',x']\in\iota_!X$, with $x'\in X(c')$ and $w'\in W_{c'}$. An arrow $V(c,x)\to s$ is an element $m\in M$ with $m\cdot[1,x]=[w',x']$. Write $m=w''g$ with $g\colon c\to d$ and $w''\in W_d$ (Lemma~\ref{lem:fact}). By Proposition~\ref{prop:free}, $m\cdot[1,x]=[w'',X(g)x]$, and by uniqueness this equals $s$ if and only if $d=c'$, $w''=w'$ and $X(g)x=x'$. Therefore $m\mapsto g$ is a bijection
\[
\Hom_{M\ltimes\iota_!X}\big(V(c,x),\,[w',x']\big)\ \cong\ \Hom_{\int X}\big((c,x),\,(c',x')\big),
\]
with inverse $g\mapsto w'g$. It is natural in $(c,x)$: precomposing with $V(f)=f$ corresponds to precomposing with $f$, because $\iota$ is a functor. Hence $s\mapsto(c',x')$ extends uniquely to a right adjoint $R$, whose counit is the arrow corresponding to $\mathrm{id}_{(c',x')}$, namely $w'\colon[1,x']\to[w',x']$. When $w'=1$ the bijection says that $V$ is fully faithful. Finally, a functor with an adjoint induces a homotopy equivalence of classifying spaces.
\end{proof}

Concretely, $R$ sends an arrow $m\colon[w,x]\to[w',x']$ to the unique $g$ with $mw=w'g$. The root of an element changes only when a whole branch is absorbed. The adjunction is the precise form of the image used in Section~\ref{sec:universal-monoid}: the types have been abolished as explicit labels, but each element of $\iota_!X$ still knows, through $R$, which element of $X$ it grew from.

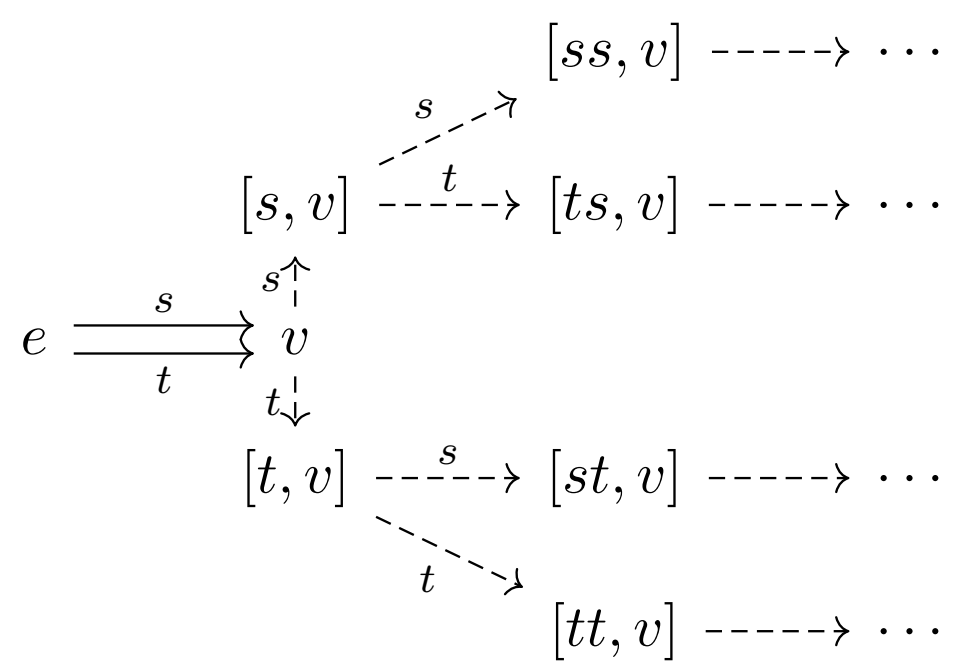
\begin{figure}[ht]
\[
\begin{tikzcd}[row sep=small, column sep=normal]
 & & {[ss,v]} \arrow[r,dashed] & \cdots\\
 & {[s,v]} \arrow[ur,dashed,"s"] \arrow[r,dashed,"t"] & {[ts,v]} \arrow[r,dashed] & \cdots\\
e \arrow[r, shift left, "s"] \arrow[r, shift right, "t"'] & v \arrow[u,dashed,"s"] \arrow[d,dashed,"t"'] & & \\
 & {[t,v]} \arrow[r,dashed,"s"] \arrow[dr,dashed,"t"'] & {[st,v]} \arrow[r,dashed] & \cdots\\
 & & {[tt,v]} \arrow[r,dashed] & \cdots
\end{tikzcd}
\]
\caption{The circle as a graph $X$ on the Kronecker category, with one edge $e$ and one vertex $v$, and the translation category of $\iota_!X$ under the free monoid $\langle s,t\rangle$. Every object emits one arrow per generator. The edge absorbs both letters (solid arrows); the vertex absorbs none, and a binary tree of branches grows from it (dashed arrows). The core $\{e\rightrightarrows v\}$ is $\int X\simeq S^1$, and the branches retract onto their roots.}
\label{fig:kron}
\end{figure}

\begin{remark}[Core and branches]\label{rem:core}
The image of $V$, the \emph{core}, is the full subcategory on the roots $[1,x]$. Every object $[w,x]$ receives exactly one arrow from the core whose image under $R$ is an identity, namely the counit $w$. If $\C$ has no retractions (no composite of non-identity arrows is an identity), no branch returns to the core: the core is a sieve, and when $\C$ is free the branches are trees, as in Figure~\ref{fig:kron}. With retractions, branches may return. If $rs=\mathrm{id}_a$ with $s\colon a\to b$, then for $x\in X(b)$ the element $[s,x]$ lies on a branch, and
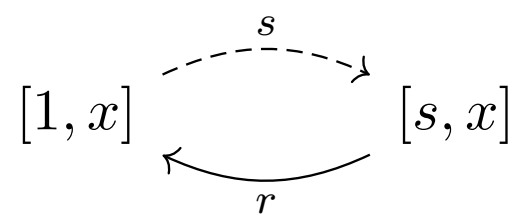
\[
\begin{tikzcd}[column sep=large]
{[1,x]} \arrow[r, bend left=25, dashed, "s"] & {[s,x]} \arrow[l, bend left=25, "r"]
\end{tikzcd}
\qquad r\cdot[s,x]=[rs,x]=[1,x].
\]
Theorem~\ref{thm:B} holds regardless: a coreflection does not need the core to be a sieve.
\end{remark}

Through Thomason's theorem, the coreflection theorem can also be read as a statement about Kan extensions.

\begin{corollary}\label{cor:kan}
For every functor $X\colon\C\to\Set$, the map $\hocolim_\C X\to\hocolim_{\Mon(\C)}\iota_!X$ induced by the functor $\iota$ and by the unit $X\to\iota^*\iota_!X$ is a homotopy equivalence.
\end{corollary}

\begin{proof}
Under Thomason's natural identification of both homotopy colimits with classifying spaces of Grothendieck constructions, this map is $B$ of the functor $\int_\C X\to\int_\C\iota^*\iota_!X\to\Mon(\C)\ltimes\iota_!X$ induced by the unit and by $\iota$. This functor sends $(c,x)$ to $[1,x]$ and $f$ to the letter $f$; it is the functor $V$ of Theorem~\ref{thm:B}, and $BV$ is a homotopy equivalence.
\end{proof}

The theorem also closes the computations left open in Section~\ref{sec:categories-topological-combinatorial}.

\begin{example}\label{ex:spectra}
The free monoids on one and two generators are universal monoids: $\mathbb N=\Mon(\{0\to1\})$ and $F_2=\Mon(\mathbf K)$, where $\mathbf K=\{E\rightrightarrows V\}$ is the Kronecker category. By Example~\ref{ex:sorts}, $\T(\mathbf K)$ consists of the homotopy types of all graphs, so Theorem~\ref{thm:B} and Example~\ref{ex:F2} give
\[
\T(F_2)=\{\text{homotopy types of graphs}\}.
\]
Figure~\ref{fig:kron} is the case of the circle. Similarly, $\T(\mathbb N)$ consists exactly of the disjoint unions of points and circles: one inclusion is Example~\ref{ex:N}, and the other follows from the $\mathbb N$-sets $\mathbb Z/k$ and $\mathbb N$ of that example, since homotopy quotients take coproducts of actions to disjoint unions. Here the coreflection theorem only yields the discrete spaces, $\T(\{0\to1\})$, and the circles have to be found by hand: the inclusion $\T(\CC)\subseteq\T(\Mon(\CC))$ may be strict.
\end{example}

The coreflection theorem only speaks about the actions of the form $\iota_!X$. It gives a lower bound for the realization spectrum of $\Mon(\CC)$, which is all we need to build master monoids, but it says nothing about how far the spectrum may grow. The last part of the example already shows that it may grow: the free monoid $\mathbb N$ realizes circles that the arrow category $\{0\to1\}$ cannot. Controlling arbitrary actions requires a different tool, akin to an amalgamation theorem for categories. We provide such a tool in the next section, and we use it for the monoids $M_n=\Mon((\Dinj^{\leq n})^{\op})$: the universal monoids of the semisimplicial category truncated at dimension $n$, generated by the face operators of dimension at most $n$. We shall see that no action of $M_n$, free or not, produces homotopy above dimension $n$.

\section{The semisimplicial master monoid and its levels}
\label{sec:semisimplicial-master-monoid}

One of the surprises of algebraic topology is that the knowledge of how to glue simplices to obtain any homotopy type is contained in the study of the injective maps between finite ordinals, which in turn admit a combinatorial presentation with very simple relations. Using this we will build a master monoid.

More precisely, in this section we apply the construction of the universal monoid to the semisimplicial category and to its truncations. The result is a sequence of monoids $M_0\to M_1\to M_2\to\cdots\to M_\Delta$ which reflects algebraically the filtration of semisimplicial sets by dimension. Each $M_n$ is finitely presented and realizes exactly the homotopy types of dimension at most $n$, while the colimit $M_\Delta$ realizes every CW homotopy type. The relation between the different dimensions recalls the theory of levels of the toposes of simplicial and cubical sets developed in \cite{KennettRiehlRoyZaks2011}, which we already met when discussing the \emph{Aufhebung}.

\subsection{The universal monoid of the semisimplicial category}
\label{subsec:semisimplicial-universal-monoid}

The category $\Dinj$ of finite nonempty ordinals $[n]=\{0<1<\cdots<n\}$, $n\geq0$, and strictly increasing maps is generated by the coface maps $\delta_i^k\colon[k-1]\to[k]$, $0\leq i\leq k$, whose image misses $i$. They satisfy the semisimplicial identities $\delta_j^k\delta_i^{k-1}=\delta_i^k\delta_{j-1}^{k-1}$ for $0\leq i<j\leq k$. We work with the opposite category $\CC_\Delta=\Dinj^{\op}$ and write $d_i^k=(\delta_i^k)^{\op}\colon[k]\to[k-1]$ for the corresponding face operators. With our convention for composition, the identities read
\[
d_i^{k-1}d_j^k=d_{j-1}^{k-1}d_i^k,\qquad 0\leq i<j\leq k.
\]
A functor $X\colon\CC_\Delta\to\Set$ is precisely a semisimplicial set.

\begin{definition}
A monoid $M$ is a \emph{master monoid} if every CW homotopy type is of the form $S/\!/M$ for some $M$-set $S$. The \emph{semisimplicial master monoid} is the universal monoid $M_\Delta=\Mon(\Dinj^{\op})$, with universal functor $\iota_\Delta\colon\Dinj^{\op}\to M_\Delta$.
\end{definition}

Applying the description of the universal monoid by generators and relations, we obtain the presentation
\[
M_\Delta=\left\langle d_i^k\ \middle|\ \begin{array}{l} k\geq1,\ 0\leq i\leq k,\\[2pt] d_i^{k-1}d_j^k=d_{j-1}^{k-1}d_i^k\quad(0\leq i<j\leq k)\end{array}\right\rangle.
\]
Thus $M_\Delta$ is obtained by identifying all the objects of $\Dinj^{\op}$ and regarding all the face operators, in all dimensions, as elements of a single monoid. The functor $\iota_{\Delta!}$ turns a semisimplicial set $X$ into the free $M_\Delta$-set $\iota_{\Delta!}X=M_\Delta\otimes_{\Dinj^{\op}}X$. Since the category of elements of $X$ has the homotopy type of its geometric realization (Example~\ref{ex:sorts}), the coreflection theorem gives a natural homotopy equivalence $|X|\simeq\iota_{\Delta!}X/\!/M_\Delta$.

\begin{theorem}[The semisimplicial master monoid]
\label{thm:semisimplicial-master}
The monoid $M_\Delta$ is a master monoid: for every CW complex $Y$ there is an $M_\Delta$-set $S$ with $S/\!/M_\Delta\simeq Y$.
\end{theorem}

\begin{proof}
Every CW complex has the homotopy type of the realization of a simplicial complex $K$ \cite[Thm.~2, case $n=1$]{Milnor59}. The barycentric subdivision $\sd K$ can be read as a semisimplicial set $X$, whose simplices are the chains of simplices of $K$ ordered by inclusion, and $|X|\cong|\sd K|\cong|K|\simeq Y$. Taking $S=\iota_{\Delta!}X$, Theorem~\ref{thm:B} gives $S/\!/M_\Delta\simeq B\!\int_{\Dinj^{\op}}X\simeq|X|\simeq Y$.
\end{proof}

The construction is strictly algebraic: the homotopy type of $Y$ is encoded by a set and an action $M_\Delta\times S\to S$. The original dimension of each element of $X$ disappears as an explicit type, but it remains recorded in the face operators that the element can absorb.

\subsection{Truncations and finite presentations}
\label{subsec:semisimplicial-truncations}

For $n\geq0$ let $\CC_n=(\Dinj^{\leq n})^{\op}$ and let $j_n\colon\CC_n\hookrightarrow\CC_\Delta$ be the inclusion. The functors $X\colon\CC_n\to\Set$ are the semisimplicial sets truncated at dimension $n$.

\begin{definition}
The \emph{semisimplicial monoid truncated at dimension $n$} is $M_n=\Mon(\CC_n)$, with universal functor $\iota_n\colon\CC_n\to M_n$.
\end{definition}

The functor $\Mon$ turns $j_n$ into a homomorphism $\mu_n=\Mon(j_n)$, and we obtain a commutative square
\[
\begin{tikzcd}[column sep=large]
\CC_n \ar[r,hook,"j_n"] \ar[d,"\iota_n"'] & \CC_\Delta \ar[d,"\iota_\Delta"]\\
M_n \ar[r,"\mu_n"'] & M_\Delta.
\end{tikzcd}
\]
More generally, for $n\leq m$ we obtain homomorphisms $\mu_{n,m}\colon M_n\to M_m$, compatible with the inclusions of the truncated categories, and a sequence $M_0\to M_1\to M_2\to\cdots$ whose colimit is $M_\Delta\cong\colim_nM_n$.

\begin{proposition}[Presentation of the truncated monoids]
\label{prop:presentation-Mn}
For every $n\geq1$, the monoid $M_n$ admits the finite presentation
\[
M_n=\left\langle d_i^k\ \middle|\ \begin{array}{l} 1\leq k\leq n,\ 0\leq i\leq k,\\[2pt] d_i^{k-1}d_j^k=d_{j-1}^{k-1}d_i^k,\\[2pt] 2\leq k\leq n,\ 0\leq i<j\leq k\end{array}\right\rangle.
\]
In particular, $M_n$ has $\sum_{k=1}^n(k+1)=n(n+3)/2$ generators and $\sum_{k=2}^n\binom{k+1}{2}$ relations in this presentation.
\end{proposition}

\begin{proof}
The category $\CC_n$ is presented by the face operators $d_i^k$, $1\leq k\leq n$, and the semisimplicial identities between parallel paths. Forming the universal monoid identifies all the objects but keeps the same generators and relations. There are $k+1$ generators in dimension $k$, and the relations involving dimensions $k$, $k-1$ and $k-2$ are indexed by the pairs $0\leq i<j\leq k$.
\end{proof}

The first examples are $M_0=\{1\}$ and $M_1=\langle d_0^1,d_1^1\rangle$, the free monoid on two generators. In dimension two we obtain
\[
M_2=\left\langle d_0^1,d_1^1,d_0^2,d_1^2,d_2^2\ \middle|\ d_0^1d_1^2=d_0^1d_0^2,\ \ d_0^1d_2^2=d_1^1d_0^2,\ \ d_1^1d_2^2=d_1^1d_1^2\right\rangle,
\]
with five generators and three relations. The rewriting system of Lemma~\ref{lem:nf} gives effective normal forms for the elements of $M_n$; since $\CC_n$ is finite, this system is finite and convergent, and in particular the word problem of $M_n$ is decidable.

\subsection{Restriction and extension between levels}
\label{subsec:restriction-between-levels}

The inclusion $j_n\colon\CC_n\hookrightarrow\CC_\Delta$ induces the usual triple of adjoints $j_{n!}\dashv j_n^*\dashv j_{n*}$. Here $j_n^*$ is the truncation functor, which keeps only the simplices of dimension at most $n$, and its adjoints are the Kan extensions $j_{n!}=\Lan_{j_n}$ and $j_{n*}=\Ran_{j_n}$. In the semisimplicial case, $j_{n!}$ extends an $n$-truncated semisimplicial set by the empty set in all dimensions above $n$, so its essential image consists of the semisimplicial sets of dimension at most $n$; on the other hand, $j_{n*}$ completes a truncated object by adding, in higher dimensions, all compatible families of faces, much as the clique complex of a graph fills in every complete subgraph.

Similarly, the homomorphism $\mu_n$ induces a triple $\mu_{n!}\dashv\mu_n^*\dashv\mu_{n*}$. The restriction functors commute strictly, $j_n^*\iota_\Delta^*=\iota_n^*\mu_n^*$, since both composites restrict an action of $M_\Delta$ to the face operators of dimension at most $n$. Taking left adjoints we obtain a canonical natural isomorphism $\mu_{n!}\iota_{n!}\cong\iota_{\Delta!}j_{n!}$, so that the square
\[
\begin{tikzcd}[column sep=large]
\Set^{\CC_n} \ar[r,"j_{n!}"] \ar[d,"\iota_{n!}"'] & \Set^{\CC_\Delta} \ar[d,"\iota_{\Delta!}"]\\
\Set^{M_n} \ar[r,"\mu_{n!}"'] & \Set^{M_\Delta}
\end{tikzcd}
\]
commutes up to isomorphism. We may therefore first form the free action of $M_n$ and then extend it to $M_\Delta$, or first extend the truncated semisimplicial set and then form the free action of $M_\Delta$: the results are canonically isomorphic. Taking right adjoints we likewise obtain $\mu_{n*}\iota_{n*}\cong\iota_{\Delta*}j_{n*}$. Thus the initial square controls simultaneously the three levels of the Kan calculus: restriction, free extension and cofree extension.

\begin{remark}[Relation with the levels of a topos]
\label{rem:levels-topos}
For a full inclusion of small categories, the triple formed by restriction and the two Kan extensions describes an essential subtopos of the presheaf category. In the simplicial, cubical and globular toposes these levels coincide with the dimensions \cite{KennettRiehlRoyZaks2011}, and the tower $\Set^{\CC_0}\to\Set^{\CC_1}\to\Set^{\CC_2}\to\cdots$ can be read in this way: the dimension $n$ determines a level of the topos of semisimplicial sets. The monoidal tower $\Set^{M_0}\to\Set^{M_1}\to\Set^{M_2}\to\cdots$ is compatible with this filtration through the squares above, but we do not claim that $\Set^{M_n}$ is a level of $\Set^{M_\Delta}$ in the technical sense. Although $j_n$ is fully faithful, the homomorphism $\mu_n$ arises after all the objects have been identified, and its Kan extensions need not be fully faithful. What survives is the formal compatibility between truncations and free actions.
\end{remark}

\subsection{Dimensional universality}
\label{subsec:dimensional-universality}

The tower of monoids $M_n$ is not only an algebraic reflection of the semisimplicial truncations. Each level controls exactly one dimensional bound in homotopy.

\begin{definition}
Let $n\geq0$. A homotopy type has \emph{dimension at most $n$}, written $\operatorname{hdim}(Y)\leq n$, if it can be represented by a CW complex of dimension at most $n$. A monoid $M$ is \emph{$n$-universal} if every homotopy type of dimension at most $n$ is of the form $S/\!/M$ for some $M$-set $S$.
\end{definition}

\begin{lemma}
\label{lem:semisimplicial-models}
A homotopy type has dimension at most $n$ if and only if it can be represented by the geometric realization $|Z|$ of a semisimplicial set $Z$ of dimension at most $n$.
\end{lemma}

\begin{proof}
The realization of a semisimplicial set of dimension at most $n$ is a CW complex of dimension at most $n$, with one $k$-cell for each $k$-simplex. Conversely, every CW complex of dimension at most $n$ is homotopy equivalent to a simplicial complex $K$ of dimension at most $n$ \cite[Thm.~2C.5]{Hatcher}; for countable complexes see also \cite[Thm.~13]{Whitehead}.
The idea is to proceed by induction on skeleta, replacing the attaching maps of the cells by homotopic simplicial maps. Finally, the barycentric subdivision of $K$, whose $k$-simplices are the chains $\sigma_0\subsetneq\cdots\subsetneq\sigma_k$ of simplices of $K$, is a semisimplicial set of dimension at most $n$ whose realization is homeomorphic to $|K|$.
\end{proof}

\begin{theorem}[Dimensional universality]
\label{thm:Mn-n-universal}
For every $n\geq1$, the monoid $M_n$ is $n$-universal.
\end{theorem}

\begin{proof}
Let $Y$ be a homotopy type of dimension at most $n$. By Lemma~\ref{lem:semisimplicial-models}, it is represented by $|Z|$ for a semisimplicial set $Z$ of dimension at most $n$, that is, with $Z_k=\varnothing$ for $k>n$. Since $Z$ is empty above dimension $n$, it is determined by its restriction $X\colon\CC_n\to\Set$, and the categories of elements $\int_{\CC_n}X$ and $\int_{\Dinj^{\op}}Z$ are isomorphic; hence $B\!\int_{\CC_n}X\simeq|Z|\simeq Y$ by Example~\ref{ex:sorts}. Applying the coreflection theorem to $\iota_n\colon\CC_n\to M_n$, the $M_n$-set $S=\iota_{n!}X$ satisfies $S/\!/M_n=B(M_n\ltimes\iota_{n!}X)\simeq B\!\int_{\CC_n}X\simeq Y$.
\end{proof}

The case $n=1$ recovers the free monoid on two generators, $M_1\cong F_2$, whose actions realize all the homotopy types of graphs (Example~\ref{ex:spectra}). In dimension two, the finitely presented monoid $M_2$, with five generators and three relations, realizes every homotopy type of dimension at most two, in particular that of every two-dimensional CW complex. Let us see the most familiar of them.

\begin{example}[The sphere from two discs]
\label{ex:sphere}
Let $X$ be the $2$-dimensional semisimplicial set with one vertex $v$, one edge $e$ with $d_0e=d_1e=v$, and two triangles $\sigma,\tau$ all of whose faces are $e$. Its realization is the circle $|e|$ with two discs attached along the same map, which runs once around the circle (replacing this attaching map by a homotopic homeomorphism does not change the homotopy type \cite[Prop.~0.18]{Hatcher}). Two discs glued along their common boundary form a sphere, so $|X|\simeq S^2$. As an $M_2$-set, $S=\iota_{2!}X$ is generated by $\sigma$ and $\tau$ subject to the relations
\[
d^2_0\sigma=d^2_1\sigma=d^2_2\sigma=d^2_0\tau=d^2_1\tau=d^2_2\tau,
\]
which identify all six faces with the edge $e$. The relation $d^1_0e=d^1_1e$ need not be imposed, since it follows from the semisimplicial identity $d^1_0d^2_2=d^1_1d^2_0$ applied to $\sigma$.
The diagram shows the core of $M_2\ltimes S$ together with every generator arrow that leaves each of its objects. Solid arrows stay in the core, and dashed arrows go to branches, which we cut after one step.
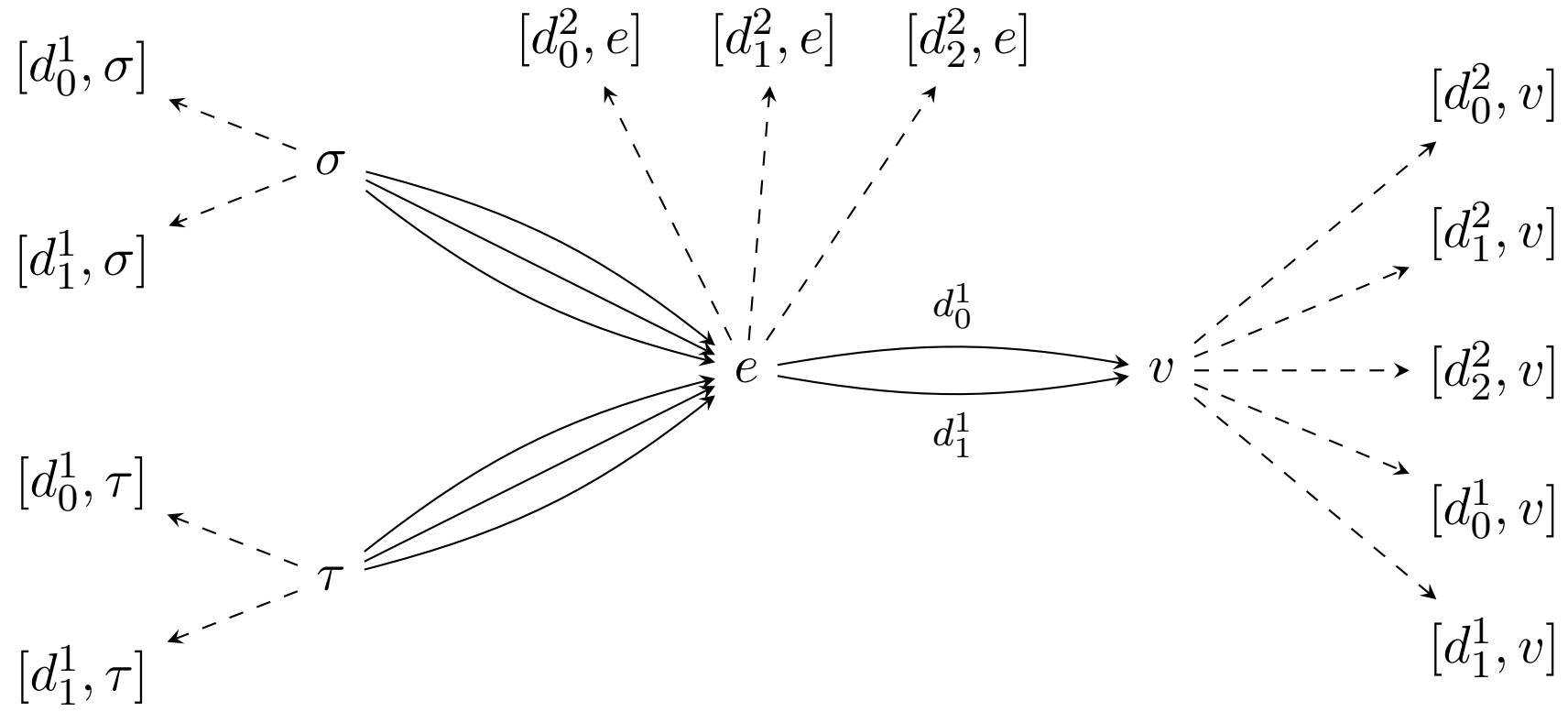
\[
\begin{tikzpicture}[x=1cm,y=1cm,>=stealth,baseline=(e.base)]
\node (s) at (0,1.5) {$\sigma$};
\node (t) at (0,-1.5) {$\tau$};
\node (e) at (3,0) {$e$};
\node (v) at (6,0) {$v$};
\node (s0) at (-1.8,2.2) {$[d^1_0,\sigma]$};
\node (s1) at (-1.8,0.8) {$[d^1_1,\sigma]$};
\node (t0) at (-1.8,-0.8) {$[d^1_0,\tau]$};
\node (t1) at (-1.8,-2.2) {$[d^1_1,\tau]$};
\node (e0) at (1.8,2.4) {$[d^2_0,e]$};
\node (e1) at (3.2,2.4) {$[d^2_1,e]$};
\node (e2) at (4.6,2.4) {$[d^2_2,e]$};
\node (v0) at (8.4,2) {$[d^2_0,v]$};
\node (v1) at (8.4,1) {$[d^2_1,v]$};
\node (v2) at (8.4,0) {$[d^2_2,v]$};
\node (v3) at (8.4,-1) {$[d^1_0,v]$};
\node (v4) at (8.4,-2) {$[d^1_1,v]$};
\draw[->] (s) to[bend left=12] (e);
\draw[->] (s) -- (e);
\draw[->] (s) to[bend right=12] (e);
\draw[->] (t) to[bend left=12] (e);
\draw[->] (t) -- (e);
\draw[->] (t) to[bend right=12] (e);
\draw[->] (e) to[bend left=10] node[above]{\scriptsize$d^1_0$} (v);
\draw[->] (e) to[bend right=10] node[below]{\scriptsize$d^1_1$} (v);
\foreach \a/\b in {s/s0,s/s1,t/t0,t/t1,e/e0,e/e1,e/e2,v/v0,v/v1,v/v2,v/v3,v/v4}
  \draw[->,dashed] (\a) -- (\b);
\end{tikzpicture}
\]
The triangles absorb the three letters $d^2_j$ and send $d^1_0,d^1_1$ to branches; the edge absorbs the letters $d^1_i$ and sends the $d^2_j$ to branches; the vertex absorbs nothing and emits five branches. Along the branches the letters still interact through the semisimplicial identities, for instance $d^1_0\cdot[d^2_0,e]=d^1_0\cdot[d^2_1,e]$. The solid part is $\int X$, and Theorem~\ref{thm:B} gives $S/\!/M_2\simeq|X|\simeq S^2$. Two triangles whose faces have been identified, acted on by a monoid with five generators and three relations, are a sphere.
\end{example}

The compatibility between levels allows us to state the construction stably. If $X$ is a semisimplicial set of dimension at most $n$, then $\mu_{n!}\iota_{n!}X\cong\iota_{\Delta!}j_{n!}X$, so the realization of $X$ through $M_n$ and its realization through $M_\Delta$ are related by extension of scalars along $\mu_n$, and $\iota_{n!}X/\!/M_n\simeq|X|\simeq\iota_{\Delta!}j_{n!}X/\!/M_\Delta$. The monoid $M_\Delta$ realizes all the levels at once, but every finite-dimensional homotopy type already appears in some truncated monoid $M_n$.

\subsection{Exactness of the truncated master monoids}
\label{subsec:exactness-truncated-monoids}

The universality above only provides an inclusion: every homotopy type of dimension at most $n$ can be realized by an action of $M_n$. In our case the converse also holds: an arbitrary action of $M_n$ cannot produce homotopy of dimension higher than $n$. Here the free actions are not enough, and we need a description of the translation categories of arbitrary actions of a universal monoid. Its proof rests on an amalgamation theorem: a statement of this kind for semigroupoids already appears in McDuff's work, as a step in the proof of her realization theorem \cite[Lemma~4]{McDuff}, and we use it in the form given by Hurwitz \cite[Thm.~2]{Hurwitz}. There is a certain justice in this: the tool that allowed McDuff to realize each space by a monoid is the same one that allows us to control all the actions of a single monoid.

\begin{lemma}\label{lem:arbitrary-actions}
Let $\C$ be a small category and let $S$ be a $\Mon(\C)$-set. Then
\[
S/\!/\Mon(\C)\simeq\hocolim\big(B\!\textstyle\int_\C\iota^*S\longleftarrow\Ob\C\times S\longrightarrow S\big),
\]
where $\Ob\C\times S$ and $S$ are discrete, the left map is the inclusion of objects and the right one is the projection. In particular, $B\Mon(\C)\simeq\hocofib(\Ob\C\to B\C)$.
\end{lemma}

\begin{proof}
Write $M=\Mon(\C)$ and $O=\Ob\C\times S$.

\emph{Step 1: $M\ltimes S$ is a pushout in $\Cat$.} The projection $P\colon M\ltimes S\to M$ is a discrete opfibration: for every object $s$ and every arrow $m$ of $M$ there is exactly one arrow of $M\ltimes S$ with source $s$ over $m$, namely $m\colon s\to m\cdot s$. In particular $P$ is a Conduché functor: if an arrow $h\colon s\to t$ lies over a composite $m_2m_1$, then $h$ factors as $s\xrightarrow{m_1}m_1\cdot s\xrightarrow{m_2}t$, and this is the only factorization of $h$ over $(m_1,m_2)$. By the theorem of Giraud and Conduché \cite{Giraud,Conduche}, the pullback functor $P^*\colon\Cat/M\to\Cat/(M\ltimes S)$ has a right adjoint, and therefore preserves colimits; colimits in a slice of $\Cat$ are computed in $\Cat$. Regard the pushout square defining $\Mon(\C)$ as a square over $M$ and apply $P^*$. Its corners become: $P^*(M)=M\ltimes S$; $P^*(*)$, the fibre of $P$ over the unique object and its identity, which is the discrete category $S$; $P^*(\Ob\C)$, the discrete category $O$; and $P^*(\C)=\C\times_M(M\ltimes S)$, whose objects are the pairs $(c,s)$ and whose arrows are the $f\colon(c,s)\to(c',f\cdot s)$ with $f\colon c\to c'$ in $\C$, that is, $\int_\C\iota^*S$. The induced functors are the inclusion of objects $O\to\int\iota^*S$, the projection $O\to S$, the functor $\int\iota^*S\to M\ltimes S$ given by $(c,s)\mapsto s$ and $f\mapsto f$, and the inclusion of objects $S\to M\ltimes S$. Hence
\[
M\ltimes S\;\cong\;\textstyle\int_\C\iota^*S\sqcup_{O}S.\tag{$*$}
\]

\emph{Step 2: a groupoid replacing $S$.} Let $\mathrm{Ch}$ be the groupoid with set of objects $O$ having exactly one arrow $(c,s)\to(c',s')$ when $s=s'$ and none otherwise. The functor $q\colon\mathrm{Ch}\to S$, $(c,s)\mapsto s$, is fully faithful and surjective on objects, hence an equivalence of categories, and $Bq\colon B\mathrm{Ch}\to S$ is a homotopy equivalence, each component of $\mathrm{Ch}$ being an indiscrete category. The projection $O\to S$ of $(*)$ is the composite of the inclusion of objects $O\to\mathrm{Ch}$ with $q$.

\emph{Step 3: the amalgam.} Let $\Cat_O$ be the category of small categories with set of objects $O$ and functors that are the identity on objects; its initial object is the discrete category $O$. Both $\int\iota^*S$ and $\mathrm{Ch}$ are objects of $\Cat_O$, and we let $E=\int\iota^*S*_O\mathrm{Ch}$ be their pushout in $\Cat_O$; since the two functors out of $O$ are bijective on objects, $E$ is also their pushout in $\Cat$. In Hurwitz's terminology \cite[\S2]{Hurwitz}, an inclusion $L\hookrightarrow N$ in $\Cat_O$ is right-free when the arrows of $N$ decompose into disjoint cosets under precomposition with arrows of $L$, with a fixed set of representatives; when $L=O$ is discrete, every arrow of $N$ forms its own coset, so the inclusions $O\hookrightarrow\int\iota^*S$ and $O\hookrightarrow\mathrm{Ch}$ are right-free. Hurwitz's amalgamation theorem \cite[Thm.~2]{Hurwitz} then shows that the canonical map
\[
B\!\textstyle\int\iota^*S\cup_{O}B\mathrm{Ch}\longrightarrow BE
\]
is a homotopy equivalence, where the union is taken along the two inclusions of $O=BO$ as sets of vertices. These are inclusions of subcomplexes, so the union is a homotopy pushout, and since $Bq$ is a homotopy equivalence under $O$, we get
\[
BE\simeq\hocolim\big(B\!\textstyle\int\iota^*S\longleftarrow O\longrightarrow S\big).
\]

\emph{Step 4: from $E$ to $M\ltimes S$.} Let $j\colon\mathrm{Ch}\to E$ be the canonical functor and form the pushout of $q$ along $j$:
\[
\begin{tikzcd}
\mathrm{Ch} \arrow[r,"j"] \arrow[d,"q"'] & E \arrow[d,"q'"]\\
S \arrow[r] & E\sqcup_{\mathrm{Ch}}S.
\end{tikzcd}
\]
By the pasting of pushouts, $E\sqcup_{\mathrm{Ch}}S\cong\int\iota^*S\sqcup_O S$, which is $M\ltimes S$ by $(*)$. We use the folk model structure on $\Cat$, whose weak equivalences are the equivalences of categories, whose cofibrations are the functors injective on objects, and whose fibrations are the isofibrations \cite{JT,Rezk}. Here $q$ is a weak equivalence and $j$, being bijective on objects, is a cofibration. Every small category is cofibrant in this model structure, since the functor from the empty category is injective on objects, and a model category in which every object is cofibrant is left proper \cite[Prop.~13.1.2]{Hirschhorn}. Hence $q'\colon E\to M\ltimes S$ is an equivalence of categories. An equivalence of categories induces a homotopy equivalence of classifying spaces, because natural isomorphisms induce homotopies; thus the weak equivalences of the folk structure are in particular functors inducing homotopy equivalences of nerves, which is all we use, and the Thomason model structure plays no role. Therefore $S/\!/M=B(M\ltimes S)\simeq BE$, which proves the first statement. For $S=*$ we have $M\ltimes*=M$ and $\int\iota^*(*)=\C$, and the homotopy pushout of $B\C\leftarrow\Ob\C\to*$ is the homotopy cofibre of the inclusion of the objects as vertices.
\end{proof}

The lemma is, in a sense, dual to the coreflection theorem. There we started from a $\C$-set $X$, whose elements carry a type, and asked how it looks when it is regarded as an action of $\Mon(\C)$; the answer was the free action $\iota_!X$, in which the types survive only as the letters that each element can absorb. Here we go in the opposite direction. We start from an arbitrary $\Mon(\C)$-set $S$, whose elements have no type at all, and look at it through the $\C$-set $\iota^*S$. Restriction does not choose a type for each element: it gives every element all the types at once, so that each $s\in S$ appears once for every object $c$ of $\C$, as the object $(c,s)$ of $\int_\C\iota^*S$, and these typed copies of $s$ are in general different objects. To recover the action we have to undo this multiplication of copies, that is, to erase the distinction between the objects of $\C$. The pushout $(*)$ does this strictly, and a pushout of categories need not have any homotopical meaning; the content of the lemma is that this one does. Geometrically, $S/\!/\Mon(\C)$ is obtained from $B\!\int\iota^*S$ by attaching, for each $s\in S$, a cone on the points $(c,s)$, $c\in\Ob\C$: a new point joined by an arc to each typed copy of $s$. The copies are not glued together; they are joined by paths, so the distinction between the objects is erased only up to homotopy. This is why forgetting the types costs, homotopically, only $1$-cells.

\begin{definition}
A monoid $M$ is \emph{$n$-exact} if it is $n$-universal and, for every $M$-set $S$, the homotopy quotient $S/\!/M$ has dimension at most $n$.
\end{definition}

\begin{theorem}[Exactness of the truncated levels]
\label{thm:Mn-n-exact}
For every $n\geq1$, the monoid $M_n$ is $n$-exact. Consequently, its homotopy quotients are exactly the homotopy types of dimension at most $n$:
\[
\T(M_n)=\left\{[Y]\ \middle|\ \operatorname{hdim}(Y)\leq n\right\}.
\]
\end{theorem}

\begin{proof}
Universality is Theorem~\ref{thm:Mn-n-universal}; it remains to prove the upper bound. Let $S$ be an $M_n$-set, let $Z=\iota_n^*S$, regarded as a semisimplicial set of dimension at most $n$, and let $O=\Ob(\CC_n)\times S$. By Lemma~\ref{lem:arbitrary-actions},
\[
S/\!/M_n\simeq\hocolim\left(B\!\int_{\CC_n}Z\longleftarrow O\longrightarrow S\right),
\]
and $B\!\int_{\CC_n}Z\simeq|Z|$ by Example~\ref{ex:sorts}. The object $([k],s)$ of $\int Z$ receives an arrow from $([0],v(s))$, where $v(s)=Z(\epsilon)(s)$ is the $0$-th vertex of the simplex $s$, $\epsilon\colon[k]\to[0]$ being the arrow of $\CC_n$ opposite to the inclusion of $0$ in $[k]$; and by the naturality of the equivalence of Example~\ref{ex:sorts}, applied to the map $\Delta^0\to Z$ classifying $v(s)$, the object $([0],v(s))$ corresponds to the vertex $v(s)$ of $|Z|$ up to homotopy. Since a homotopy pushout of maps out of a discrete space only depends on these maps up to homotopy, $S/\!/M_n$ is homotopy equivalent to the double mapping cylinder of $O\to|Z|$, $([k],s)\mapsto v(s)$, and of the projection $O\to S$. This double mapping cylinder is the realization of the semisimplicial set obtained from $Z$ by adding a new vertex $\bar s$ for each $s\in S$ and, for each $([k],s)\in O$, an edge from $v(s)$ to $\bar s$. Its dimension is at most $\max\{n,1\}=n$, so its realization is a CW complex of dimension at most $n$, and hence $\operatorname{hdim}(S/\!/M_n)\leq n$.
\end{proof}

\begin{corollary}
\label{cor:Mn-not-higher-master}
For every $n\geq1$, the monoid $M_n$ is not $(n+1)$-universal.
\end{corollary}

\begin{proof}
Every homotopy quotient $S/\!/M_n$ has dimension at most $n$. The sphere $S^{n+1}$ does not. Hence $M_n$ cannot realize $S^{n+1}$.
\end{proof}

Exactness allows us to read the tower $M_1\to M_2\to\cdots\to M_\Delta$ as an algebraic filtration of homotopy by dimension. The level $M_n$ contains exactly the realizing power needed for the $n$-dimensional homotopy types, while the colimit $M_\Delta$ gathers all the levels:
\[
\T(M_1)\subseteq\T(M_2)\subseteq\cdots\subseteq\T(M_\Delta)=\{\text{all CW homotopy types}\}.
\]

\begin{remark}[Finite acyclic categories]
\label{rem:acyclic}
The construction applies directly to finite acyclic categories, for which we refer the reader to \cite[Chapter~10]{Kozlov2008}. Recall that a category $\CC$ is \emph{acyclic} if its only endomorphisms and its only isomorphisms are the identities. The chains $c_0\xrightarrow{f_1}c_1\to\cdots\xrightarrow{f_k}c_k$ of non-identity morphisms of $\CC$ form a semisimplicial set $N_+\CC$ (a trisp, in the terminology of \cite[Chapter~10]{Kozlov2008}), whose face operators drop one end of a chain or compose two consecutive morphisms; acyclicity guarantees that a composite of non-identity morphisms is never an identity. The category of elements of $N_+\CC$ is the subdivision $\sd\CC$, and $B(\sd\CC)\simeq B\CC$ \cite{delHoyo2008}. If $\CC$ has no chain of more than $n$ composable non-identity morphisms, $N_+\CC$ has dimension at most $n$, and the coreflection theorem gives
\[
B\CC\simeq\iota_{n!}N_+\CC/\!/M_n .
\]
Thus every finite acyclic category is realized by a finitely generated action of the finitely presented monoid $M_n$: the types of the chains are identified, and the classifying space does not change.
\end{remark}

\section{From test categories to master monoids: realizing but not modelling}
\label{sec:test-categories-master-monoids}

Whoever looks at current homotopy theory will notice that the idea of a topological space, as it is presented in courses on general topology, rarely appears explicitly. Homotopy has more to do with combinatorial models, and it is often expressed in the language of model categories rather than in that of spaces. One of the most ambitious attempts to reformulate our way of understanding homotopy types is due to Grothendieck in \emph{Pursuing Stacks} \cite{PS}, and it is known as the theory of test categories (see \cite{Malt,Cis}, and \cite{zawadowski_kapulkin_test_categories} for an accessible introduction). It inverts the usual perspective. Instead of taking a space as the fundamental object and then looking for a combinatorial description of it, one starts with a small category $\A$ whose objects play the role of elementary shapes, such as simplices, cubes, globes or trees. A presheaf on $\A$ is a space built with these shapes, its category of elements gathers all the occurrences of the shapes, and its classifying space assembles them homotopically. A test category is a category of shapes rich enough for this language to be complete.

This allows us to distinguish two levels of universality. The first is a property of representation: every homotopy type is realized by some presheaf on $\A$. The second is a property of modelling: the category of presheaves, localized at the morphisms that induce homotopy equivalences, is equivalent to the homotopy category of spaces. The second statement contains the first, but it also demands that presheaves and their morphisms reproduce homotopy theory correctly, and not only its collection of objects. The universal monoid keeps the first level, but not the second. This is a general phenomenon: the coreflection theorem shows that every small category $\CC$ transfers to $\Mon(\CC)$ its ability to realize homotopy types, $\T(\CC)\subseteq\T(\Mon(\CC))$, so that, if $\A$ is a weak test category, $M_\A=\Mon(\A^{\op})$ realizes every homotopy type through its actions on sets, and the semisimplicial monoid of the previous section is not an isolated case. But the passage to the monoid identifies all the shapes of $\A$ with a single object, and the category of all $M_\A$-sets contains many actions that do not arise from presheaves; their homotopical behaviour is no longer governed by the test property of $\A$. We shall arrive at a fundamental distinction: a monoid can \emph{realize} every homotopy type through its actions without \emph{modelling} homotopy theory.

\subsection{Test categories and homotopy types}
\label{subsec:test-categories-homotopy-types}

Let $\A$ be a small category and $\widehat{\A}=\Set^{\A^{\op}}$ its category of presheaves; we identify each object $a$ of $\A$ with the representable presheaf $\A(-,a)$, and write $e$ for the terminal presheaf. For a presheaf $X$, the objects of its category of elements $\A/X$ are the pairs $(a,x)$ with $x\in X(a)$, and a morphism $(a,x)\to(b,y)$ is a morphism $f\colon a\to b$ of $\A$ with $X(f)(y)=x$. With our conventions for covariant categories of elements there is a natural isomorphism $\A/X\cong(\int_{\A^{\op}}X)^{\op}$, and since a category and its opposite have canonically homeomorphic classifying spaces, $B(\A/X)\cong B\!\int_{\A^{\op}}X$. The functor $i_\A\colon\widehat{\A}\to\Cat$, $X\mapsto\A/X$, has a right adjoint $i_\A^*$, which sends a small category $\mathcal B$ to the presheaf $a\mapsto\Hom_{\Cat}(\A/a,\mathcal B)$ \cite[1.2.1 and 1.2.2]{Malt}. We follow Maltsiniotis \cite{Malt}, and take as definitions some of the equivalent characterizations given there.

\begin{definition}
\label{def:test}
Let $\A$ be a small category.
\begin{enumerate}[leftmargin=*]
\item A morphism $u\colon X\to Y$ of presheaves is a \emph{weak equivalence} if $B(\A/X)\to B(\A/Y)$ is a homotopy equivalence. We denote by $W_\A$ the class of weak equivalences.
\item The category $\A$ is \emph{aspherical} if $B\A\simeq*$. A presheaf $X$ is \emph{aspherical} if $\A/X$ is aspherical, and \emph{locally aspherical} if the presheaf $a\times X$ is aspherical for every object $a$ of $\A$ \cite[1.2.6]{Malt}.
\item $\A$ is a \emph{weak test category} if the presheaf $i_\A^*(\mathcal B)$ is aspherical for every small category $\mathcal B$ with a terminal object \cite[Prop.~1.3.9]{Malt}.
\item $\A$ is a \emph{local test category} if the slice $\A/a$ is a weak test category for every object $a$, and a \emph{test category} if it is both a weak test category and a local test category \cite[Def.~1.5.2]{Malt}.
\item $\A$ is \emph{totally aspherical} if it is aspherical and the presheaf $a\times b$ is aspherical for all objects $a,b$ \cite[Def.~1.6.2]{Malt}. A \emph{strict test category} is a totally aspherical test category \cite[Def.~1.6.7]{Malt}.
\end{enumerate}
\end{definition}

If $\A$ is a weak test category, the functor $X\mapsto\A/X$ induces an equivalence $W_\A^{-1}\widehat{\A}\simeq\Hot$, where $\Hot$ is the homotopy category of spaces \cite[\S1.3]{Malt}; this is the precise sense in which presheaves on $\A$ model homotopy types. A weak test category is aspherical \cite[Prop.~1.3.9]{Malt}. This can be understood by looking at the terminal presheaf: its category of elements is $\A/e\cong\A$, so for the terminal object of $\widehat{\A}$ to represent a point it is necessary that $B\A\simeq*$.

\begin{remark}
Test category theory usually works with presheaves, that is, with contravariant functors $X\colon\A^{\op}\to\Set$. The reason seems to us rather historical. The founding example is the category $\Delta$ of simplices, which models spaces precisely when we take the opposite category, since the face maps between topological simplices go in the opposite direction to the face operators of a simplicial set. Moreover, sheaves and presheaves were built with the poset of open sets in mind, where an inclusion $U\hookrightarrow V$ becomes a restriction $X(V)\to X(U)$; and the contravariant form is the usual one when working with the Yoneda embedding.
\end{remark}

The test condition is hard to check directly. Grothendieck found a criterion that puts its whole weight on a single object.

\begin{definition}[{\cite[1.4.1]{Malt}}]\label{def:interval}
An \emph{interval} of $\A$ is a presheaf $I$ together with two global sections $\partial_0,\partial_1\colon e\to I$. It is \emph{separating} if the equalizer of $\partial_0$ and $\partial_1$ is the empty presheaf.
\end{definition}

\begin{theorem}[Grothendieck]
\label{thm:interval}
Let $\A$ be a small category.
\begin{enumerate}[leftmargin=*]
\item $\A$ is a local test category if and only if it admits a separating interval $I$ which is locally aspherical; equivalently, if and only if its Lawvere object, the subobject classifier of $\widehat{\A}$, is locally aspherical \cite[Thm.~1.5.6]{Malt}.
\item $\A$ is a test category if and only if it is aspherical and admits a separating, locally aspherical interval \cite[Rem.~1.5.4]{Malt}.
\item If $\A$ is totally aspherical, then $\A$ is a weak test category if and only if it is a strict test category, and this happens if and only if $\A$ admits a separating interval $I$ which is aspherical \cite[Prop.~1.6.8]{Malt}.
\end{enumerate}
\end{theorem}

It is worth pausing on what is asked of the interval. It plays the role of the unit interval $[0,1]$, that is, of the cylinder of the point: for every presheaf $X$, the product $X\times I$ is a cylinder of $X$, with its two ends $X\rightrightarrows X\times I$ given by $\partial_0$ and $\partial_1$, and a homotopy is a map out of $X\times I$. Each of the two conditions protects one half of this role. Local asphericity means precisely that every projection $X\times I\to X$ is a weak equivalence \cite[1.2.6]{Malt}: crossing with $I$ does not change the homotopy type, so that the interval is, homotopically, a point. Separation means that the two ends are distinct in the strongest possible sense: there is no shape $a$ of $\A$ at which $\partial_0$ and $\partial_1$ agree. One may read this, in the dialectical vocabulary of Lawvere that we already met when speaking of the \emph{Aufhebung} \cite{LawvereFutureCategoryTheory}, as a unity and distinction of opposites. The two ends are opposites which are never identified, not even locally, and yet they are held together in a single object which is homotopically a point. Neither half suffices on its own: the point $e$ with $\partial_0=\partial_1$ is united but not distinct, and the two points $e\sqcup e$ are distinct but not united, since $a\times(e\sqcup e)=a\sqcup a$ has a disconnected category of elements. Homotopy lives in the tension between both: a path is what makes two different things the same up to homotopy without identifying them. This is, in a way, the opposite of what we did with the universal monoid, where the objects of a category were identified; we shall see in Section~\ref{subsec:test-monoid} that a monoid with a good interval, in which the two constants $0$ and $1$ play the role of the opposite ends, does model homotopy theory.

\begin{example}
The simplicial category $\Delta$ is the fundamental example of a test category. Joyal's category $\Theta$ is also a test category, as shown by Cisinski and Maltsiniotis \cite{CM}. These cases provide different systems of combinatorial shapes for representing homotopy types through presheaves. The semisimplicial category $\Dinj$ used in the previous section is a weak test category \cite[Prop.~1.7.25]{Malt}, although it is not a test category (see the introduction of \cite{Malt}).
\end{example}

In this paper we mainly use the essential consequence of the weak test condition: every homotopy type is represented by the category of elements of some presheaf on $\A$. Hence, if $\A$ is a weak test category, $\T(\A^{\op})$ consists of all homotopy types. This equality only expresses the essential surjectivity on the objects of $\Hot$; the definition of a weak test category contains a stronger statement, since it also demands that the morphisms of the homotopy category be recovered through the localization of $\widehat{\A}$ at $W_\A$.

\subsection{From test categories to master monoids}
\label{subsec:from-test-to-master}

Let now $\A$ be a weak test category. We apply the construction of the universal monoid to the opposite category $\CC=\A^{\op}$, and write $M_\A=\Mon(\A^{\op})$, with universal functor $\iota_\A\colon\A^{\op}\to M_\A$. Each presheaf $X$ is an $\A^{\op}$-set, and its free extension along $\iota_\A$ is the $M_\A$-set $\iota_{\A!}X=M_\A\otimes_{\A^{\op}}X$.

\begin{theorem}[From test categories to master monoids]
\label{thm:test-category-master-monoid}
If $\A$ is a weak test category, then $M_\A=\Mon(\A^{\op})$ is a master monoid.
\end{theorem}

\begin{proof}
Let $Y$ be a homotopy type. Since $\A$ is a weak test category, there is a presheaf $X$ with $B(\A/X)\simeq Y$. Applying the coreflection theorem to $\iota_\A$ we get $B(\A/X)\cong B\!\int_{\A^{\op}}X\simeq B(M_\A\ltimes\iota_{\A!}X)$, so the $M_\A$-set $S=\iota_{\A!}X$ satisfies $S/\!/M_\A\simeq Y$.
\end{proof}

\begin{corollary}
\label{cor:test-master-examples}
The monoids $\Mon(\Delta^{\op})$ and $\Mon(\Theta^{\op})$ are master monoids. Moreover, $M_\Delta=\Mon(\Dinj^{\op})$ is a master monoid by Theorem~\ref{thm:semisimplicial-master}; this also follows from Theorem~\ref{thm:test-category-master-monoid}, since $\Dinj$ is a weak test category \cite[Prop.~1.7.25]{Malt}.
\end{corollary}

In particular, $M_\Delta=\Mon(\Dinj^{\op})$ is a master monoid, as we already proved directly using semisimplicial sets and barycentric subdivisions. Theorem~\ref{thm:test-category-master-monoid} shows that this phenomenon does not depend on the particular properties of simplices. It is a formal consequence of two facts:
\begin{enumerate}
\item the presheaves on a weak test category realize every homotopy type;
\item the coreflection theorem transfers each realization from a small category to its universal monoid.
\end{enumerate}

\subsection{Objectwise universality}
\label{subsec:objectwise-universality}

It is worth making precise which part of the test condition has been used in the proof above. We have only used that $X\mapsto B(\A/X)$ is essentially surjective onto homotopy types, and not the stronger statement $W_\A^{-1}\widehat{\A}\simeq\Hot$. This leads to an intermediate notion.

\begin{definition}
A small category $\CC$ is \emph{objectwise universal} if $\T(\CC)$ consists of all homotopy types, that is, if for every homotopy type $Y$ there is a functor $X\colon\CC\to\Set$ with $B\!\int_{\CC}X\simeq Y$.
\end{definition}

If $\A$ is a weak test category, then $\A^{\op}$ is objectwise universal, but the converse need not hold. Objectwise universality only says that all the objects of $\Hot$ appear as realizations; it provides no description of the morphisms between them, nor does it claim that a localization of $\Set^{\CC}$ is equivalent to $\Hot$.

\begin{proposition}
\label{prop:Mon-preserves-objectwise-universality}
If a small category $\CC$ is objectwise universal, then so is $\Mon(\CC)$.
\end{proposition}

\begin{proof}
By the coreflection theorem, $\T(\CC)\subseteq\T(\Mon(\CC))$.
\end{proof}

The coreflection theorem gives more than an abstract inclusion between classes. It defines an explicit transformation $X\mapsto\iota_!X$, and the original category of elements appears as a full coreflective subcategory of the translation category, so that the original categorical model can be recovered canonically inside the free action. For a weak test category $\A$ we can summarize the situation by the diagram
\[
\begin{tikzcd}[column sep=large]
\Set^{\A^{\op}} \ar[r,"\iota_{\A!}"] \ar[d,"X\mapsto B(\A/X)"'] & \Set^{M_\A} \ar[d,"S\mapsto S/\!/M_\A"]\\
\Hot \ar[r,equal] & \Hot,
\end{tikzcd}
\]
which commutes up to natural homotopy equivalence. The top functor transfers representatives of homotopy types, but this commutativity does not imply that it induces an equivalence between the localizations of the two categories of actions.

\subsection{The failure of asphericity}
\label{subsec:failure-asphericity}

The first obstruction for the universal monoid to be a weak test category appears when we look at its own classifying space. Even if $\A$ is aspherical, the monoid $\Mon(\A^{\op})$ need not be. Identifying the objects may create new cohomology classes in dimension one.

\begin{proposition}
\label{prop:universal-monoid-not-aspherical}
Let $\CC$ be a connected small category with at least two objects. Then $H^1(B\Mon(\CC);\mathbb Z)\neq0$. In particular, $\Mon(\CC)$ is not aspherical, and it is not a weak test category.
\end{proposition}

\begin{proof}
By Lemma~\ref{lem:arbitrary-actions}, $B\Mon(\CC)$ is the homotopy cofibre of the inclusion $\Ob\CC\to B\CC$ of the objects as vertices. Its long exact sequence in reduced homology contains
\[
H_1(B\CC)\longrightarrow H_1(B\Mon(\CC))\longrightarrow\tilde H_0(\Ob\CC)\longrightarrow\tilde H_0(B\CC).
\]
Since $\CC$ is connected, $\tilde H_0(B\CC)=0$, so $H_1(B\Mon(\CC))$ maps onto $\tilde H_0(\Ob\CC)$, the kernel of the augmentation $\mathbb Z^{(\Ob\CC)}\to\mathbb Z$, which is a nonzero free abelian group because $\CC$ has at least two objects. Hence $H^1(B\Mon(\CC);\mathbb Z)\cong\Hom(H_1(B\Mon(\CC)),\mathbb Z)\neq0$, by the universal coefficient theorem. Therefore $B\Mon(\CC)$ is not contractible, and since every weak test category is aspherical \cite[Prop.~1.3.9]{Malt}, $\Mon(\CC)$ is not a weak test category.
\end{proof}

More precisely, the homotopy cofibre is $B\CC$ with a new point $p$ and an arc from $p$ to each vertex $c$. Choosing an object $c_0$ and, for every other object $c$, a path in $B\CC$ from $c$ to $c_0$, and sliding the end at $c$ of each arc along it \cite[Prop.~0.18]{Hatcher}, we get a noncanonical homotopy equivalence
\[
B\Mon(\CC)\simeq B\CC\vee\bigvee_{c\in\Ob\CC\setminus\{c_0\}}S^1 .
\]
Identifying all the objects of $\CC$ adds one circle for each object but one: these circles record the loss of information about the types, and they are the paths along which Lemma~\ref{lem:arbitrary-actions} erases, up to homotopy, the distinction between the objects.

\begin{corollary}
\label{cor:master-not-test}
Let $\A$ be a weak test category with at least two objects. Then $\Mon(\A^{\op})$ is a master monoid which is not a weak test category.
\end{corollary}

\begin{proof}
Theorem~\ref{thm:test-category-master-monoid} shows that $\Mon(\A^{\op})$ is a master monoid. On the other hand, $\A$ is aspherical, hence connected, and Proposition~\ref{prop:universal-monoid-not-aspherical} applied to $\CC=\A^{\op}$ shows that $\Mon(\A^{\op})$ is not aspherical, hence not a weak test category.
\end{proof}

\begin{corollary}
\label{cor:semisimplicial-master-not-test}
The semisimplicial master monoid $M_\Delta=\Mon(\Dinj^{\op})$ realizes every homotopy type, but it is not a weak test category: $\Dinj^{\op}$ is connected and has infinitely many objects, so $H^1(BM_\Delta;\mathbb Z)\neq0$ by Proposition~\ref{prop:universal-monoid-not-aspherical}.
\end{corollary}

The difference between realizing and modelling can be stated precisely. For a monoid $M$, consider the functor $\mathcal R_M\colon\Set^M\to\Hot$, $S\mapsto S/\!/M$. If $M$ is a master monoid, $\mathcal R_M$ is essentially surjective: every homotopy type is isomorphic in $\Hot$ to $\mathcal R_M(S)$ for some $M$-set $S$. Essential surjectivity, however, does not imply that $\mathcal R_M$ presents $\Hot$ as a localization of $\Set^M$. For the one-object category $M$ to be a weak test category it would be necessary, in particular, that $BM\simeq*$, and the universal monoids attached to test categories do not satisfy this condition\footnote{Strictly speaking, the presheaves on $M$ are the right $M$-sets, while our $M$-sets are left actions, that is, presheaves on $M^{\op}$; since $BM\cong BM^{\op}$ and $\Mon(\CC)^{\op}\cong\Mon(\CC^{\op})$, the obstruction does not depend on this choice.}. The situation can be summarized as follows: if $\A$ is a weak test category, then $\Mon(\A^{\op})$ realizes every homotopy type, while in general $\Mon(\A^{\op})$ is not a weak test category and the failure can already be detected at the terminal object.

\begin{remark}[Where the point lives]
\label{rem:point}
The terminal $M$-set, the one-point set with the trivial action, does not represent a point, since $*/\!/M=BM$. The point is nevertheless realized, by the free action on one generator, that is, by the regular $M$-set $M$. In its translation category the unit $1$ is an initial object: from this single root grows one branch for each generator of $M$, then, from every element of a branch, one arrow for each generator again, and so on; the relations of $M$ only make some of these branches meet. Hence $M/\!/M\simeq*$. For a universal monoid this is the free action generated by a representable functor, $\iota_!\CC(c,-)\cong M$, by the co-Yoneda lemma; for $M_\Delta$ it is the free action generated by a vertex, or by any simplex. The diagram shows the beginning of $M_2\ltimes M_2$, where the semisimplicial identity $d^1_0d^2_1=d^1_0d^2_0$ makes two branches meet:
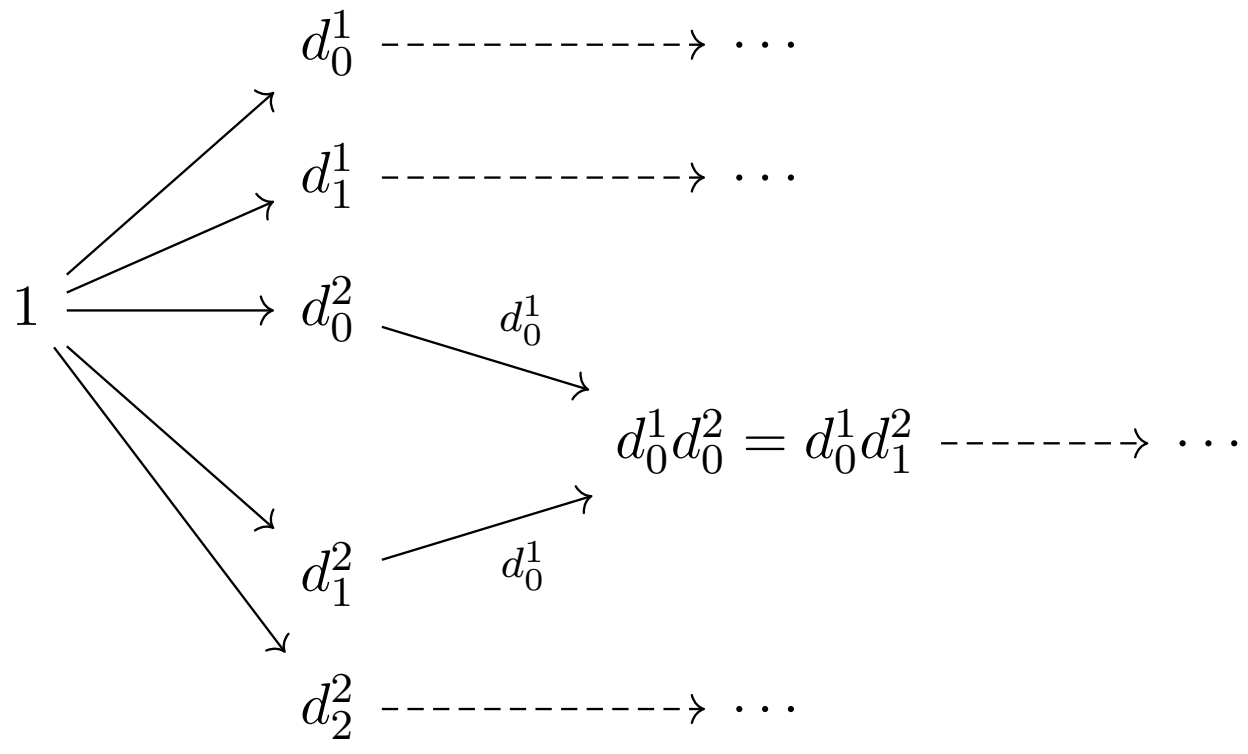
\[
\begin{tikzcd}[row sep=tiny, column sep=large]
 & d^1_0 \arrow[r,dashed] & \cdots \\
 & d^1_1 \arrow[r,dashed] & \cdots \\
1 \arrow[uur] \arrow[ur] \arrow[r] \arrow[ddr] \arrow[dddr] & d^2_0 \arrow[dr,"d^1_0"] & \\
 & & d^1_0d^2_0=d^1_0d^2_1 \arrow[r,dashed] & \cdots \\
 & d^2_1 \arrow[ur,"d^1_0"'] & \\
 & d^2_2 \arrow[r,dashed] & \cdots
\end{tikzcd}
\]
If $\A$ is a weak test category, the point is also realized by the free action generated by the terminal presheaf $1$: by the coreflection theorem, $\iota_{\A!}1/\!/M_\A\simeq B\!\int_{\A^{\op}}1\cong B\A^{\op}\simeq*$. Thus in $\Set^{M_\A}$ the point is represented by free actions, while the terminal object, which in a test category plays the role of the point, represents $BM_\A$.
\end{remark}

\subsection{Splitting instead of identifying}
\label{subsec:test-monoid}

The obstructions above might suggest that a single monoid can realize homotopy types but never model them. This is not so. What fails in Corollary~\ref{cor:master-not-test} is not the passage to a single object, but the way in which it is made. When we identify the objects of a test category, the shapes are abolished, and the circles of $B\Mon(\A^{\op})$ record their disappearance.

There is another way of reducing a category to a monoid, in which the shapes are not abolished but retained as idempotents. If every object of a category $\mathcal C$ is a retract of a single object $g$, then the category of presheaves on $\mathcal C$ is equivalent to the category of presheaves on the monoid $\operatorname{End}_{\mathcal C}(g)$, since each object is recovered as the image of an idempotent endomorphism of $g$.\footnote{Equivalently, $\mathcal C$ and the one-object category $\operatorname{End}_{\mathcal C}(g)$ have equivalent Karoubi envelopes. Since categories of $\mathbf{Set}$-valued presheaves are invariant under Karoubi completion, their presheaf categories are equivalent (see \cite[Theorem~1]{BorceuxDejean1986}).}

The classical example is Lawvere's: reflexive graphs are the right actions of a monoid with three elements, the identity and the two idempotents that send an edge to its source and to its target \cite{LawvereGraphs}. Along this second path, a single monoid can also model homotopy theory.

Recall that presheaves on a monoid $E$, regarded as a one-object category, are right $E$-sets; with our convention for actions they are the $E^{\op}$-sets, and $E$ is a weak test category exactly when $\mathcal R_{E^{\op}}\colon\Set^{E^{\op}}\to\Hot$ presents $\Hot$ as a localization.

\begin{proposition}\label{prop:test-monoid}
Let $X$ be an infinite set and let $E=\mathrm{End}(X)$ be the monoid of all maps $X\to X$, with product $gf=g\circ f$. Then $E$, regarded as a one-object category, is a strict test category. In particular, $E^{\op}$ is a master monoid whose actions also model homotopy theory.
\end{proposition}

\begin{proof}
We work in ZFC: the axiom of choice is needed only to identify $X\times X$ and $X\times\{0,1\}$ with $X$, which is where the hypothesis that $X$ is infinite enters. We check the criteria of Theorem~\ref{thm:interval}, with the notions of Definition~\ref{def:test}: $E$ is aspherical, it has a separating interval which is locally aspherical, and it is totally aspherical. Recall that presheaves on $E$ are right $E$-sets, that the terminal presheaf is the one-point set, that the only representable presheaf is $E$ with right multiplication, and that a presheaf $F$ is aspherical when its category of elements $E/F$ is.

\emph{$E$ is aspherical.} Let $K\colon E\to E$ be the functor sending every element to $1$; it factors through the trivial monoid. An element $z\in E$, regarded as an endomorphism of the unique object, is a natural transformation $\mathrm{Id}_E\Rightarrow K$ if and only if, for every $f\in E$, the naturality square commutes, that is, $K(f)\circ z=z\circ\mathrm{Id}_E(f)$, or $z=z\circ f$ with our convention $gf=g\circ f$. Every constant map $z$ satisfies $z\circ f=z$ for all $f$ (whereas $f\circ z$ is constant but in general different from $z$). Hence a constant map gives a natural transformation $\mathrm{Id}_E\Rightarrow K$, so the identity of $BE$ is homotopic to a constant map, and $BE$ is contractible.

\emph{A separating interval.} Let $I$ be the set of maps $\varphi\colon X\to\{0,1\}$, with the right action $\varphi\cdot f=\varphi\circ f$. The two constant maps are fixed by $E$, so they define global sections $\partial_0,\partial_1\colon*\to I$; since they are different, their equalizer is the empty presheaf, and $(I,\partial_0,\partial_1)$ is a separating interval (Definition~\ref{def:interval}).

\emph{Products.} Since $X$ is infinite, the axiom of choice provides bijections $\beta\colon X\times X\to X$ and $\gamma\colon X\times\{0,1\}\to X$. The maps
\[
E\times E\to E,\ (a,b)\mapsto\beta\circ(a,b),\qquad E\times I\to E,\ (a,\varphi)\mapsto\gamma\circ(a,\varphi),
\]
are bijections, with inverses
\[
m\mapsto(\mathrm{pr}_1\beta^{-1}m,\ \mathrm{pr}_2\beta^{-1}m),\qquad m\mapsto(\mathrm{pr}_1\gamma^{-1}m,\ \mathrm{pr}_2\gamma^{-1}m),
\]
and they are isomorphisms of right $E$-sets, for the diagonal actions on the products, because $(a,b)\circ f=(a\circ f,b\circ f)$ gives $\beta\circ(a\circ f,b\circ f)=(\beta\circ(a,b))\circ f$, and similarly for $\gamma$. 

The representable presheaf is aspherical, since its category of elements has a terminal object. Since $E$ is the only representable presheaf, $E\times E$ is the only product of two representables, and the first isomorphism makes it representable, hence aspherical; together with the asphericity of $E$, this means that $E$ is totally aspherical (Definition~\ref{def:test}). Likewise $E\times I\cong E$ is aspherical, and since $E$ is the only representable presheaf, this is exactly the condition that $I$ be locally aspherical. By Theorem~\ref{thm:interval}(2), $E$ is a test category, and being totally aspherical, it is a strict test category (Definition~\ref{def:test}).
\end{proof}

The proof uses only the constant maps, the maps with values in $\{0,1\}$, and the bijections $\beta,\gamma$ together with the projections of their inverses. This gives a countable example.

\begin{corollary}\label{cor:test-monoid-pr}
Let $P$ be the monoid of primitive recursive functions $\mathbb N\to\mathbb N$ under composition. Then $P$, regarded as a one-object category, is a strict test category.
\end{corollary}

\begin{proof}
The set $P$ is closed under composition and contains the constant maps, so the asphericity argument of Proposition~\ref{prop:test-monoid} applies. Let $I_P\subseteq P$ be the right $P$-set of primitive recursive functions with values in $\{0,1\}$, with the action by precomposition; the constants $0$ and $1$ are global sections, and since they are different their equalizer is empty, so $I_P$ is a separating interval. We use two different bijections: Cantor's pairing function and a second, simpler one,
\[
\begin{aligned}
\beta&\colon\mathbb N\times\mathbb N\to\mathbb N, & \beta(a,b)&=\tfrac12(a+b)(a+b+1)+b,\\
\gamma&\colon\mathbb N\times\{0,1\}\to\mathbb N, & \gamma(a,i)&=2a+i.
\end{aligned}
\]
Both are primitive recursive, being built from sums, products and halving of an even number. The inverse of $\gamma$ is $m\mapsto(\lfloor m/2\rfloor,\,m\bmod 2)$, with primitive recursive components. For $\beta$, write $T(k)=k(k+1)/2$ and let
\[
w(z)=\sum_{k=1}^{z}\chi\big(T(k)\leq z\big),\qquad \beta^{-1}(z)=\big(w(z)-z+T(w(z)),\ z-T(w(z))\big),
\]
where $\chi$ is the characteristic function of the relation $\leq$. Bounded sums of primitive recursive functions are primitive recursive, so both components of $\beta^{-1}$ are. The formula is correct: if $z=\beta(a,b)$ and $w=a+b$, then $T(w)\leq z\leq T(w)+w<T(w+1)$, so $w$ is the number of $k\geq1$ with $T(k)\leq z$, all of which satisfy $k\leq T(k)\leq z$; that is, $w=w(z)$, and then $b=z-T(w)$ and $a=w-b$.

Now, for $a,b\in P$ and $\varphi\in I_P$, the functions $\beta\circ(a,b)$ and $\gamma\circ(a,\varphi)$ are composites of primitive recursive functions, hence belong to $P$; for $m\in P$, the two components of $\beta^{-1}\circ m$ belong to $P$, the first component of $\gamma^{-1}\circ m$ belongs to $P$, and the second, $m\bmod 2$, takes values in $\{0,1\}$ and belongs to $I_P$. The maps $P\times P\to P$, $(a,b)\mapsto\beta\circ(a,b)$, and $P\times I_P\to P$, $(a,\varphi)\mapsto\gamma\circ(a,\varphi)$, are therefore bijections, with inverses given by the components of $\beta^{-1}\circ m$ and $\gamma^{-1}\circ m$, and they are equivariant for precomposition exactly as in the proof of Proposition~\ref{prop:test-monoid}. Since $P$ is the only representable presheaf, the remaining steps of that proof apply: $P$ is aspherical and totally aspherical, and $I_P$ is a separating, locally aspherical interval.
\end{proof}

\begin{remark}
Behind the proposition there is a test category with many objects. A small skeleton of the category of nonempty sets of cardinality at most that of $X$, with all maps, has finite products and a representable separating interval $\{0,1\}$, so it is a strict test category \cite[Ex.~1.6.4 and Cor.~1.6.9]{Malt}; every object is a retract of $X$, and its presheaves are the right $E$-sets. The shapes, the finite sets $\{0,\dots,n\}$, have not been identified: they survive inside $E$ as the idempotent retractions of $X$ onto its finite subsets. Here we have not used the universal monoid, but the simpler construction which takes the monoid of endomorphisms of a single object.
\end{remark}

We do not know whether there is a finitely generated, a finitely presented, or a finite monoid which is a test category.

\subsection{Universality without testability}
\label{subsec:universality-without-testability}

The results above separate two different levels of universality:
\begin{enumerate}
\item \emph{Universality of realization:} every homotopy type appears as $S/\!/M$ for some action of $M$ on a set.
\item \emph{Universality as a homotopical model:} the localization of $\Set^M$ at the equivalences detected by the categories of elements is equivalent to $\Hot$.
\end{enumerate}
The second property implies the first, but the converse is false. The monoids $\Mon(\A^{\op})$ attached to weak test categories provide a natural family of counterexamples: they realize every homotopy type, but when $\A$ has more than one object they are not even aspherical.

The positive relation can be expressed through the composite $\Set^{\A^{\op}}\xrightarrow{\iota_{\A!}}\Set^{M_\A}\xrightarrow{\mathcal R_{M_\A}}\Hot$. For each presheaf $X$, the coreflection theorem provides a natural equivalence $\mathcal R_{M_\A}(\iota_{\A!}X)\simeq B(\A/X)$, so that $\iota_{\A!}$ transfers to the actions of $M_\A$ all the homotopical representatives provided by the presheaves on $\A$. However, $\Set^{M_\A}$ contains many actions outside the essential image of $\iota_{\A!}$, and the homotopical behaviour of these actions is not controlled by the test theory of $\A$.

The cause of this difference lies in the forgetting of types. In a presheaf on $\A$, each element keeps the shape of $\A$ to which it belongs. After applying $\iota_{\A!}$, all the elements become part of a single set, and the operations that could not act on an element generate new branches. For the free actions coming from presheaves, the coreflection theorem shows that these branches are homotopically invisible: the category of elements remains as a full and coreflective core inside the translation category. For an arbitrary action of the monoid, on the contrary, there need not be a typed core onto which to coreflect it.

The passage $\A^{\op}\mapsto\Mon(\A^{\op})$ therefore keeps the ability to realize homotopy types, but in general it does not keep the test property. In other words, it keeps the objects of homotopy theory, in the sense that all of them can be represented, but not necessarily the relations between them that allow us to rebuild homotopy theory as a localization. To realize all spaces is not yet to model homotopy. This distinction also clarifies the sense of the term \emph{Aufhebung} used above. When all the objects of a category are identified, the types are eliminated as explicit information, but they do not disappear completely. In the free actions they remain inscribed in the operations that each element can absorb, and they can be recovered through the coreflection. The construction does not simply forget the types: it turns them from an external classification of the elements into an internal structure of their action.

From this point of view, the master monoid is not a new category of shapes replacing simplices, cubes or the objects of a test category. It is rather the result of compressing all those shapes into a single object, while keeping enough operations to rebuild, action by action, any homotopy type. The price of this compression is that the global organization of the shapes is no longer available for arbitrary actions. This price, however, is the price of identification and not of unity: when the shapes are kept as idempotents instead of being identified, a single monoid can both realize and model homotopy types (Proposition~\ref{prop:test-monoid}).

The initial question about what topology can do for algebra thus receives a provisional answer, or at least a possible line of development. Topology allows us to endow the actions of a monoid with an interesting structure. We not only have a way of attaching a topological space to an algebraic structure; we also know how to attach a topological space to each action of a monoid. Our contribution has been to study which homotopy types appear, and how, in some cases, one can check that every homotopy type appears in this way.

In turn, this sheds some light on topological spaces and on homotopy theory, which can be seen from a new perspective. Homotopy, at the purely objectual level, is the theory of how a specific monoid, such as any monoid attached to a test category, acts on sets. Topology is the dynamics of transformations, when these transformations are those produced by the way objects are glued according to a test category. Where the test category distinguished shapes, the monoid distinguishes possibilities of action. For example, dimension no longer appears as a label attached to each element; it remains latent in the operations that the element can absorb.

Shape has not simply been forgotten: in the new language, the study of dimensions and elementary figures is transmuted into the study of how the different pieces can be transformed and assembled. The monoid remains one, but its actions are many. It is large; it contains multitudes of homotopy types.

\bibliographystyle{amsplain}
\bibliography{master-monoid-arXiv}

@book{Poly,
  author    = {Ara, D. and Burroni, A. and Guiraud, Y. and Malbos, P. and M{\'e}tayer, F. and Mimram, S.},
  title     = {Polygraphs: from rewriting to higher categories},
  series    = {London Math. Soc. Lecture Note Ser.},
  volume    = {495},
  publisher = {Cambridge University Press},
  address   = {Cambridge},
  year      = {2025}
}

@article{BDH,
  author  = {Baumslag, G. and Dyer, E. and Heller, A.},
  title   = {The topology of discrete groups},
  journal = {J. Pure Appl. Algebra},
  volume  = {16},
  year    = {1980},
  pages   = {1--47}
}

@book{Cis,
  author    = {Cisinski, D.-C.},
  title     = {Les pr{\'e}faisceaux comme mod{\`e}les des types d'homotopie},
  series    = {Ast{\'e}risque},
  volume    = {308},
  publisher = {Soci{\'e}t{\'e} Math{\'e}matique de France},
  address   = {Paris},
  year      = {2006}
}

@article{CM,
  author  = {Cisinski, D.-C. and Maltsiniotis, G.},
  title   = {La cat{\'e}gorie {$\Theta$} de {J}oyal est une cat{\'e}gorie test},
  journal = {J. Pure Appl. Algebra},
  volume  = {215},
  year    = {2011},
  pages   = {962--982}
}

@article{Conduche,
  author  = {Conduch{\'e}, F.},
  title   = {Au sujet de l'existence d'adjoints {\`a} droite aux foncteurs ``image r{\'e}ciproque'' dans la cat{\'e}gorie des cat{\'e}gories},
  journal = {C. R. Acad. Sci. Paris S{\'e}r. A-B},
  volume  = {275},
  year    = {1972},
  pages   = {A891--A894}
}

@article{EM45,
  author  = {Eilenberg, S. and Mac~Lane, S.},
  title   = {General theory of natural equivalences},
  journal = {Trans. Amer. Math. Soc.},
  volume  = {58},
  year    = {1945},
  pages   = {231--294}
}

@article{Fied,
  author  = {Fiedorowicz, Z.},
  title   = {Classifying spaces of topological monoids and categories},
  journal = {Amer. J. Math.},
  volume  = {106},
  year    = {1984},
  pages   = {301--350}
}

@book{SevenSketches,
  author    = {Fong, B. and Spivak, D. I.},
  title     = {An invitation to applied category theory: seven sketches in compositionality},
  publisher = {Cambridge University Press},
  address   = {Cambridge},
  year      = {2019}
}

@book{Giraud,
  author    = {Giraud, J.},
  title     = {M{\'e}thode de la descente},
  series    = {M{\'e}m. Soc. Math. France},
  volume    = {2},
  publisher = {Soci{\'e}t{\'e} Math{\'e}matique de France},
  address   = {Paris},
  year      = {1964}
}

@book{PS,
  author    = {Grothendieck, A.},
  title     = {Pursuing stacks ({\`A} la poursuite des champs). {V}olume {I}},
  series    = {Doc. Math. (Paris)},
  volume    = {20},
  publisher = {Soci{\'e}t{\'e} Math{\'e}matique de France},
  address   = {Paris},
  year      = {2022},
  note      = {edited by G. Maltsiniotis; manuscript written in 1983}
}

@book{Hatcher,
  author    = {Hatcher, A.},
  title     = {Algebraic topology},
  publisher = {Cambridge University Press},
  address   = {Cambridge},
  year      = {2002}
}

@book{Hirschhorn,
  author    = {Hirschhorn, P. S.},
  title     = {Model categories and their localizations},
  series    = {Math. Surveys Monogr.},
  volume    = {99},
  publisher = {Amer. Math. Soc.},
  address   = {Providence, RI},
  year      = {2003}
}

@article{Milnor59,
  author  = {Milnor, J.},
  title   = {On spaces having the homotopy type of a {CW}-complex},
  journal = {Trans. Amer. Math. Soc.},
  volume  = {90},
  year    = {1959},
  pages   = {272--280}
}

@article{Whitehead,
  author  = {Whitehead, J. H. C.},
  title   = {Combinatorial homotopy. {I}},
  journal = {Bull. Amer. Math. Soc.},
  volume  = {55},
  year    = {1949},
  pages   = {213--245}
}

@article{Hurwitz,
  author  = {Hurwitz, C. M.},
  title   = {On the homotopy theory of monoids},
  journal = {J. Austral. Math. Soc. Ser. A},
  volume  = {47},
  year    = {1989},
  pages   = {171--185}
}

@incollection{JT,
  author    = {Joyal, A. and Tierney, M.},
  title     = {Strong stacks and classifying spaces},
  booktitle = {Category theory ({C}omo, 1990)},
  series    = {Lecture Notes in Math.},
  volume    = {1488},
  publisher = {Springer},
  address   = {Berlin},
  year      = {1991},
  pages     = {213--236}
}

@article{KT,
  author  = {Kan, D. M. and Thurston, W. P.},
  title   = {Every connected space has the homology of a {$K(\pi,1)$}},
  journal = {Topology},
  volume  = {15},
  year    = {1976},
  pages   = {253--258}
}

@article{KennettRiehlRoyZaks2011,
  author  = {Kennett, C. and Riehl, E. and Roy, M. and Zaks, M.},
  title   = {Levels in the toposes of simplicial sets and cubical sets},
  journal = {J. Pure Appl. Algebra},
  volume  = {215},
  year    = {2011},
  pages   = {949--961}
}

@book{kromer2007tool,
  author    = {Kr{\"o}mer, R.},
  title     = {Tool and object: a history and philosophy of category theory},
  series    = {Sci. Netw. Hist. Stud.},
  volume    = {32},
  publisher = {Birkh{\"a}user},
  address   = {Basel},
  year      = {2007}
}

@incollection{LawvereFutureCategoryTheory,
  author    = {Lawvere, F. W.},
  title     = {Some thoughts on the future of category theory},
  booktitle = {Category theory ({C}omo, 1990)},
  series    = {Lecture Notes in Math.},
  volume    = {1488},
  publisher = {Springer},
  address   = {Berlin},
  year      = {1991},
  pages     = {1--13}
}

@incollection{LawvereGraphs,
  author    = {Lawvere, F. W.},
  title     = {Qualitative distinctions between some toposes of generalized graphs},
  booktitle = {Categories in computer science and logic ({B}oulder, {CO}, 1987)},
  series    = {Contemp. Math.},
  volume    = {92},
  publisher = {Amer. Math. Soc.},
  address   = {Providence, RI},
  year      = {1989},
  pages     = {261--299}
}

@article{LEE,
  author  = {Lee, M.-J.},
  title   = {Homotopy for functors},
  journal = {Proc. Amer. Math. Soc.},
  volume  = {36},
  year    = {1972},
  pages   = {571--577}
}

@article{LEE2,
  author  = {Lee, M.-J.},
  title   = {Erratum to: ``{H}omotopy for functors''},
  journal = {Proc. Amer. Math. Soc.},
  volume  = {42},
  year    = {1974},
  pages   = {648--650}
}

@book{Malt,
  author    = {Maltsiniotis, G.},
  title     = {La th{\'e}orie de l'homotopie de {G}rothendieck},
  series    = {Ast{\'e}risque},
  volume    = {301},
  publisher = {Soci{\'e}t{\'e} Math{\'e}matique de France},
  address   = {Paris},
  year      = {2005}
}

@article{McDuff,
  author  = {McDuff, D.},
  title   = {On the classifying spaces of discrete monoids},
  journal = {Topology},
  volume  = {18},
  year    = {1979},
  pages   = {313--320}
}

@unpublished{Rezk,
  author = {Rezk, C.},
  title  = {A model category for categories},
  note   = {Preprint},
  year   = {2000}
}

@book{RiehlCategoryTheoryContext,
  author    = {Riehl, E.},
  title     = {Category theory in context},
  publisher = {Dover Publications},
  address   = {Mineola, NY},
  year      = {2016}
}

@article{Schreier,
  author  = {Schreier, O.},
  title   = {Die {U}ntergruppen der freien {G}ruppen},
  journal = {Abh. Math. Sem. Univ. Hamburg},
  volume  = {5},
  year    = {1927},
  pages   = {161--183}
}

@article{SpivakDatabaseQueries,
  author  = {Spivak, D. I.},
  title   = {Database queries and constraints via lifting problems},
  journal = {Math. Structures Comput. Sci.},
  volume  = {24},
  year    = {2014},
  pages   = {Art. e240602, 55 pp.}
}

@article{SpivakDataMigration,
  author  = {Spivak, D. I.},
  title   = {Functorial data migration},
  journal = {Inform. and Comput.},
  volume  = {217},
  year    = {2012},
  pages   = {31--51}
}

@article{BorceuxDejean1986,
  author  = {Borceux, Francis and Dejean, Dominique},
  title   = {Cauchy completion in category theory},
  journal = {Cahiers de Topologie et G{\'e}om{\'e}trie Diff{\'e}rentielle Cat{\'e}goriques},
  volume  = {27},
  number  = {2},
  year    = {1986},
  pages   = {133--146}
}

@article{delHoyo2008,
  author  = {del Hoyo, Mat{\'\i}as L.},
  title   = {On the subdivision of small categories},
  journal = {Topology Appl.},
  volume  = {155},
  number  = {11},
  year    = {2008},
  pages   = {1189--1200},
  doi     = {10.1016/j.topol.2008.02.006}
}

@book{Kozlov2008,
  author    = {Kozlov, Dmitry},
  title     = {Combinatorial Algebraic Topology},
  series    = {Algorithms and Computation in Mathematics},
  volume    = {21},
  publisher = {Springer},
  address   = {Berlin, Heidelberg},
  year      = {2008}
}

@article{Th,
  author  = {Thomason, R. W.},
  title   = {Homotopy colimits in the category of small categories},
  journal = {Math. Proc. Cambridge Philos. Soc.},
  volume  = {85},
  year    = {1979},
  pages   = {91--109}
}

@misc{zawadowski_kapulkin_test_categories,
  author       = {Zawadowski, M.},
  title        = {Introduction to test categories},
  howpublished = {Notes by C. Kapulkin of a talk by M. Zawadowski, available at \url{https://www.math.uwo.ca/faculty/kapulkin/notes/test_categories.pdf}}
}

\end{document}